# THE LEGENDRE TRANSFORM IN MODERN OPTIMIZATION

ROMAN A. POLYAK

ABSTRACT. The Legendre transform (LET) is a product of a general duality principle: any smooth curve is, on the one hand, a locus of pairs, which satisfy the given equation and, on the other hand, an envelope of a family of its tangent lines.

An application of the LET to a strictly convex and smooth function leads to the Legendre identity (LEID). For strictly convex and three times differentiable function the LET leads to the Legendre invariant (LEINV).

Although the LET has been known for more then 200 years both the LEID and the LEINV are critical in modern optimization theory and methods.

The purpose of the paper (survey) is to show the role of the LEID and the LEINV play in both constrained and unconstrained optimization.

## 1. INTRODUCTION

Application of the duality principle to a strictly convex $f : \mathbb{R} \to \mathbb{R}$, leads to the Legendre transform

$$f^*(s) = \sup_{x \in \mathbb{R}}\{sx - f(x)\},$$

which is often called the Legendre-Fenchel transform (see [21],[29],[30]).

The LET, in turn, leads to two important notions: the Legendre identity

$$f^{*'}(s) \equiv f^{'-1}(s)$$

and the Legendre invariant

$$\text{LEINV}(f) = \left| \frac{d^3 f}{dx^3} \left( \frac{d^2 f}{dx^2} \right)^{-\frac{3}{2}} \right| = \left| -\frac{d^3 f^*}{ds^3} \left( \frac{d^2 f^*}{ds^2} \right)^{-\frac{3}{2}} \right|.$$

Our first goal is to show a number of duality results for optimization problems with equality and inequality constraints obtained in a unified manner by using LEID.

A number of methods for constrained optimization, which have been introduced in the past several decades and for a long time seemed to be unconnected, turned out to be equivalent. We start with two classical methods for equality constrained optimization.

First, the primal penalty method by Courant [16] and its dual equivalent - the regularization method by Tichonov [60].

Second, the primal multipliers method by Hestenes [28] and Powell [52], and its dual equivalent - the quadratic proximal point method by Moreau [38], Martinet [35], [36] Rockafellar [56]-[57] (see also [2], [7], [24], [27], [45] and references therein).

Classes of primal SUMT and dual interior regularization, primal nonlinear rescaling (NR) and dual proximal points with $\varphi$- divergence distance functions, primal

Lagrangian transformation (LT) and dual interior ellipsoids methods turned out to be equivalent.

We show that LEID is a universal tool for establishing the equivalence results, which are critical, for both understanding the nature of the methods and establishing their convergence properties.

Our second goal is to show how the equivalence results can be used for convergence analysis of both primal and dual methods.

In particular, the primal NR method with modified barrier (MBF) transformation leads to the dual proximal point method with Kullback-Leibler entropy divergence distance (see [50]). The corresponding dual multiplicative algorithm, which is closely related to the EM method for maximum likelihood reconstruction in position emission tomography as well as to image space reconstruction algorithm (see [17], [20], [62]), is the key instrument for establishing convergence of the MBF method (see [31], [46], [50], [53]).

In the framework of LT the MBF transformation leads to the dual interior proximal point method with Bregman distance (see [39], [49]).

The kernel $\varphi(s) = -\ln s + s - 1$ of the Bregman distance is a self-concordant (SC) function. Therefore the corresponding interior ellipsoids are Dikin's ellipsoids.

Application LT for linear programming (LP) calculations leads to Dikin's type method for the dual LP (see [18]).

The SC functions have been introduced by Yuri Nesterov and Arkadi Nemirovski in the late 80s (See [42],[43]).

Their remarkable SC theory is the centerpiece of the interior point methods (IPMs), which for a long time was the main stream in modern optimization. The SC theory establishes the IPMs complexity for large classes of convex optimization problem from a general and unique point of view.

It turns out that a strictly convex $f \in C^3$ is self-concordant if LEINV($f$) is bounded. The boundedness of LEINV($f$) leads to the basic differential inequality, four sequential integrations of which produced the main SC properties.

The properties, in particular, lead to the upper and lower bounds for $f$ at each step of a special damped Newton method for unconstrained minimization SC functions. The bounds allow establishing global convergence and show the efficiency of the damped Newton method for minimization SC function.

The critical ingredients in these developments are two special SC function: $w(t) = t - \ln(t + 1)$ and its $LET$ $w^*(s) = -s - \ln(1 - s)$.

Usually two stages of the damped Newton method is considered (see [43]). At the first stage at each step the error bound $\Delta f(x) = f(x) - f(x^*)$ is reduced by $w(\lambda)$, where $0 < \lambda < 1$ is the Newton decrement. At the second stage $\Delta f(x)$ converges to zero with quadratic rate. We considered a middle stage where $\Delta f(x)$ converges to zero with superlinear rate, which is explicitly characterized by $w(\lambda)$ and $w^*(\lambda)$.

To show the role of LET and LEINV($f$) in unconstrained optimization of SC functions was our third goal.

The paper is organized as follows.

In the next section along with LET we consider LEID and LEINV.

In section 3 penalty and multipliers methods and their dual equivalents applied for optimization problems with equality constraints.

In section 4 the classical SUMT methods and their dual equivalents - the interior regularization methods - are applied to convex optimization problem.

In section 5 we consider the Nonlinear Rescaling theory and methods, in particular, the MBF and its dual equivalent - the prox with Kullback-Leibler entropy divergence distance.

In section 6 the Lagrangian transform (LT) and its dual equivalent - the interior ellipsoids method - are considered. In particular, the LT with MBF transformation, which leads to the dual prox with Bregman distance.

In section 7 we consider LEINV, which leads to the basic differential inequality, the main properties of the SC functions and eventually to the damped Newton method.

We conclude the paper (survey) with some remarks, which emphasize the role of LET, LEID and LEINV in modern optimization.

## 2. Legendre Transformation

We consider LET for a smooth and strictly convex scalar function of a scalar argument $f : \mathbb{R} \to \mathbb{R}$.

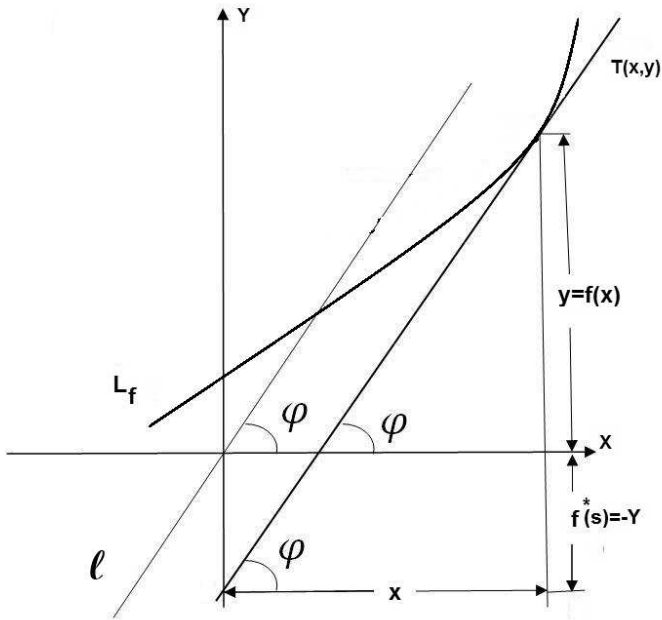

Figure 1. Legendre transformation

For a given $s = \tan \varphi$ let us consider line $l = \{(x, y) \in \mathbb{R}^2 : y = sx\}$. The corresponding tangent to the curve $L_f$ with the same slope is defined as follows:

$$T(x, y) = \{(X, Y) \in \mathbb{R}^2 : Y - f(x) = f^{'}(x)(X - x) = s(X - x)\}.$$

In other words $T(x, y)$ is a tangent to the curve $L_f = \{(x, y) : y = f(x)\}$ at the point $(x, y)$: $f'(x) = s$. For $X = 0$, we have $Y = f(x) - sx$. The conjugate function $f^* : (a, b) \to \mathbb{R}$, $-\infty < a < b < \infty$ at the point $s$ is defined as $f^*(s) = -Y = -f(x) + sx$. Therefore (see Fig. 1)

$$(1) \qquad\qquad f^*(s) + f(x) = sx.$$

More often $f^*$ is defined as follows

$$(2) \qquad\qquad f^*(s) = \max_{x \in \mathbb{R}} \{sx - f(x)\}.$$

Keeping in mind that $T(x, y)$ is the supporting hyperplane to the epi $f = \{(y, x) : y \geq f(x)\}$ the maximum in (2) is reached at $x$: $f'(x) = s$, therefore the primal representation of (1) is

$$(3) \qquad\qquad f^*(f'(x)) + f(x) \equiv f'(x)x.$$

For a strictly convex $f$ we have $f''(x) > 0$, therefore due to the Inverse Function Theorem the equation $f'(x) = s$ can be solved for $x$, that is

$$(4) \qquad\qquad x(s) = f'^{-1}(s).$$

Using (4) from (3) we obtain the dual representation of (1)

$$(5) \qquad\qquad f^*(s) + f(x(s)) \equiv sx(s).$$

Also, it follows from $f''(x) > 0$ that $x(s)$ in (2) is unique, so $f^*$ is as smooth as $f$. The variables $x$ and $s$ are not independent, they are linked through equation $s = f'(x)$.

By differentiating (5) we obtain

$$(6) \qquad\qquad f^{*\prime}(s) + f'(x(s))x'(s) \equiv x(s) + sx'(s).$$

In view of $f'(x(s)) = s$, from (4) and (6) we obtain the following identity,

$$(7) \qquad\qquad f^{*\prime}(s) \equiv f'^{-1}(s),$$

which is called the Legendre identity (LEID).

From (4) and (7) we obtain

$$(8) \qquad\qquad \frac{df^*(s)}{ds} = x.$$

On the other hand, we have

$$(9) \qquad\qquad \frac{df(x)}{dx} = s.$$

From (8) and (9) it follows

$$(10) \qquad\qquad a)\frac{d^2f^*(s)}{ds^2} = \frac{dx}{ds} \text{ and } b)\frac{d^2f(x)}{dx^2} = \frac{ds}{dx}.$$

From

$$\frac{dx}{ds} \cdot \frac{ds}{dx} = 1$$

and (10) we get

$$(11) \qquad\qquad \frac{d^2f^*}{ds^2} \cdot \frac{d^2f}{dx^2} = 1,$$

so the local curvatures of $f$ and $f^*$ are inverses to each other.

The following Theorem established the relations of the third derivatives of $f$ and $f^*$, which leads to the notion of Legendre invariant.

**Theorem 1.** *If $f \in C^3$ is strictly convex then*

$$\frac{d^3 f^*}{ds^3} \cdot \left(\frac{d^2 f^*}{ds^2}\right)^{-3/2} + \frac{d^3 f}{dx^3} \cdot \left(\frac{d^2 f}{dx^2}\right)^{-3/2} = 0. \tag{12}$$

*Proof.* By differentiating (11) in $x$ we obtain

$$\frac{d^3 f^*}{ds^3} \cdot \frac{ds}{dx} \cdot \frac{d^2 f}{dx^2} + \frac{d^2 f^*}{ds^2} \cdot \frac{d^3 f}{dx^3} = 0.$$

In view of (10b) we have

$$\frac{d^3 f^*}{ds^3} \cdot \left(\frac{d^2 f}{dx^2}\right)^2 + \frac{d^2 f^*}{ds^2} \cdot \frac{d^3 f}{dx^3} = 0. \tag{13}$$

By differentiating (11) in $s$ and keeping in mind (10a) we obtain

$$\frac{d^3 f^*}{ds^3}\frac{d^2 f}{dx^2} + \left(\frac{d^2 f^*}{ds^2}\right)^2 \frac{d^3 f}{dx^3} = 0. \tag{14}$$

Using (11), from (13) and (14) we have

$$\frac{d^3 f^*}{ds^3} \cdot \frac{d^2 f}{dx^2} + \frac{1}{\left(\frac{d^2 f}{dx^2}\right)^2}\frac{d^3 f}{dx^3} = 0$$

or

$$\frac{d^3 f^*}{ds^3} \left(\frac{d^2 f}{dx^2}\right)^3 + \frac{d^3 f}{dx^3} = 0.$$

Keeping in mind $\frac{d^2 f}{dx} > 0$ from the last equation follows

$$\frac{d^3 f^*}{ds^3} \left(\frac{d^2 f}{dx^2}\right)^{\frac{3}{2}} + \frac{d^3 f}{dx^3} \left(\frac{d^2 f}{dx^2}\right)^{-\frac{3}{2}} = 0.$$

Using (11) again we obtain (12).

**Corollary 2.** *From (12) we have*

$$-\frac{d^3 f^*}{ds^3} \left(\frac{d^2 f^*}{ds^2}\right)^{-3/2} = \frac{d^3 f}{dx^3} \left(\frac{d^2 f}{dx^2}\right)^{-3/2}.$$

*The Legendre Invariant is defined as follows*

$$\text{LEINV}(f) = \left|-\frac{d^3 f^*}{ds^3} \left(\frac{d^2 f^*}{ds^2}\right)^{-3/2}\right| = \left|\frac{d^3 f}{dx^3} \left(\frac{d^2 f}{dx^2}\right)^{-3/2}\right|. \tag{15}$$

For a strictly convex $f \in C^3$ boundedness of $\text{LEINV}(f)$ defines the class of self-concordant (SC) functions introduced by Yuri Nesterov and A. Nemirovski in the late 80s .

## 3. Equality Constrained Optimization

Let $f$ and all $c_i\colon \mathbb{R}^n \to \mathbb{R}$, $i = 1, ..., m$ be continuously differentiable. We consider the following optimization problem with equality constrains

$$
\begin{aligned}
&\min f(x) \\
&\text{s. t. } c_i(x) = 0, i = 1, ..., m.
\end{aligned}
\tag{16}
$$

We assume that (16) has a regular solution $x^*$ that is

$$
\operatorname{rank} \nabla c(x^*) = m < n,
$$

where $\nabla c(x)$ is the Jacobian of the vector - function $c(x) = (c_1(x), ..., c_m(x))^T$. Then (see, for example [45]) there exists $\lambda^* \in \mathbb{R}^m$:

$$
\nabla_x L(x^*, \lambda^*) = 0, \ \nabla_\lambda L(x^*, \lambda^*) = c(x^*) = 0,
$$

where

$$
L(x, \lambda) = f(x) + \sum_{i=1}^{m} \lambda_i c_i(x)
$$

is the classical Lagrangian, which corresponds to (16).

It is well known that the dual function

$$
d(\lambda) = \inf\{L(x, \lambda) | x \in \mathbb{R}^n\}
\tag{17}
$$

is closed and concave. Its subdifferential

$$
\partial d(\lambda) = \{g : d(u) - d(\lambda) \le (g, u - \lambda), \forall u \in \mathbb{R}^m\}
\tag{18}
$$

at each $\lambda \in \mathbb{R}^n$ is a non - empty, bounded and convex set. If for a given $\lambda \in \mathbb{R}^m$ the minimizer

$$
x(\lambda) = \operatorname{argmin}\{L(x, \lambda) | x \in \mathbb{R}^n\}
$$

exists then

$$
\nabla_x L(x(\lambda), \lambda) = 0.
\tag{19}
$$

If the minimizer $x(\lambda)$ is unique, then the dual function

$$
d(\lambda) = L(x(\lambda), \lambda)
$$

is differentiable and the dual gradient

$$
\nabla d(\lambda) = \nabla_x L(x(\lambda), \lambda) \nabla_\lambda x(\lambda) + \nabla_\lambda L(x(\lambda), \lambda),
$$

where $\nabla_\lambda x(\lambda)$ is the Jacobian of vector - function $x(\lambda) = (x_1(\lambda), ..., x_n(\lambda))^T$. In view of (19) we have

$$
\nabla d(\lambda) = \nabla_\lambda L(x(\lambda), \lambda) = c(x(\lambda)).
\tag{20}
$$

In other words, the gradient of the dual function coincides with the residual vector computed at the primal minimizer $x(\lambda)$.

If $x(\lambda)$ is not unique, then for any $\hat{x} = x(\lambda) \in \operatorname{Argmin}\{L(x, \lambda) | x \in \mathbb{R}^n\}$ we have

$$
c(\hat{x}) \in \partial d(\lambda).
$$

In fact, let

$$
u : d(u) = L(x(u), u) = \min_{x \in \mathbb{R}^n} L(x, u),
\tag{21}
$$

then for any $\lambda \in \mathbb{R}^m$ we have

$$d(u) = \min\{f(x) + \sum_{i=1}^m u_i c_i(x) | x \in \mathbb{R}^n\} \leq f(\hat{x}) + \sum_{i=1}^m u_i c_i(\hat{x}) = f(\hat{x}) + \sum \lambda_i c_i(\hat{x})$$
$$+ (c(\hat{x}), u - \lambda) = d(\lambda) + (c(\hat{x}), u - \lambda)$$

or

$$d(u) - d(\lambda) \leq (c(\hat{x}), u - \lambda), \forall u \in \mathbb{R}^m,$$

so (18) holds for $g = c(\hat{x})$, therefore

$$(22) \qquad\qquad c(\hat{x}) \in \partial d(\hat{\lambda}).$$

The dual to (16) problem is

$$(23) \qquad\qquad \max d(\lambda)$$
$$\text{s. t. } \lambda \in \mathbb{R}^m,$$

which is a convex optimization problem independent from convexity properties of $f$ and $c_i$, $i = 1, ..., m$ in (16).

The following inclusion

$$(24) \qquad\qquad 0 \in \partial d(\lambda^*)$$

is the optimality condition for the dual maximizer $\lambda^*$ in (23).

3.1. **Penalty Method and its Dual Equivalent.** In this subsection we consider two methods for solving optimization problems with equality constraints and their dual equivalents.

In 1943 Courant introduced the following penalty function and correspondent method for solving (16) (see [16]).

Let $\pi(t) = \frac{1}{2} t^2$ and $k > 0$ be the penalty (scaling) parameter, then Courant's penalty function $P : \mathbb{R}^n \times \mathbb{R}_{++} \to \mathbb{R}$ is defined by the following formula

$$(25) \qquad P(x, k) = f(x) + k^{-1} \sum_{i=1}^m \pi(k c_i(x)) = f(x) + \frac{k}{2} \|c(x)\|^2,$$

where $\|\cdot\|$ is Euclidian norm. At each step the penalty method finds unconstrained minimizer

$$(26) \qquad\qquad x(k) : P(x(k), k) = \min_{x \in \mathbb{R}^n} P(x, k).$$

We assume that for a given $k > 0$ minimizer $x(k)$ exists and can be found from the system $\nabla_x P(x, k) = 0$. Then

$$\nabla_x P(x(k), k) =$$
$$(27) \qquad \nabla f(x(k)) + \sum_{i=1}^m \pi^{'}(k c_i(x(k))) \nabla c_i(x(k)) = 0.$$

Let

$$(28) \qquad\qquad \lambda_i(k) = \pi^{'}(k c_i(x(k))), \ i = 1, .., m.$$

From (27) and (28) follows

$$(29) \quad \nabla_x P(x(k), k) = \nabla f(x(k)) + \sum_{i=1}^m \lambda_i(k) \nabla c_i(x(k)) = \nabla_x L(x(k), \lambda(k)) = 0,$$

which means that $x(k)$ satisfies the necessary condition to be a minimizer of $L(x, \lambda(k))$. If $L(x(k), \lambda(k)) = \min_{x \in \mathbb{R}^n} L(x, \lambda(k))$, then $d(\lambda(k)) = L(x(k), \lambda(k))$ and

$$(30) \qquad\qquad c(x(k)) \in \partial d(\lambda(k)).$$

Due to $\pi''(t) = 1$ the inverse function $\pi'^{-1}$ exists. From (28) follows

$$(31) \qquad\qquad c_i(x(k)) = k^{-1} \pi'^{-1}(\lambda_i(k)), \ i = 1, ..., m.$$

From (30), (31) and the LEID $\pi'^{-1} = \pi^{*'}$ we obtain

$$(32) \qquad\qquad 0 \in \partial d(\lambda(k)) - k^{-1} \sum_{i=1}^m \pi^{*'}(\lambda_i(k)) e_i,$$

where $e_i = (0, ..., 1, .., 0)$.

The inclusion (32) is the optimality condition for $\lambda(k)$ to be the unconstrained maximizer of the following unconstrained maximization problem

$$(33) \qquad d(\lambda(k)) - k^{-1} \sum_{i=1}^m \pi^*(\lambda_i(k)) = \max\{d(u) - k^{-1} \sum_{i=1}^m \pi^*(u_i) : u \in \mathbb{R}^m\}.$$

Due to $\pi^*(s) = \max_t \{st - \frac{1}{2} t^2\} = \frac{1}{2} s^2$ the problem (33) one can rewrite as follows

$$(34) \qquad d(\lambda(k)) - \frac{1}{2k} \sum_{i=1}^m \lambda_i^2(k) = \max\{d(u) - \frac{1}{2k} \|u\|^2 : u \in \mathbb{R}^m\}.$$

Thus, Courant's penalty method (26) is equivalent to Tikhonov's (see [60]) regularization method (34) for the dual problem (23).

The convergence analysis of (34) is simple because the dual $d(u)$ is concave and $D(u, k) = d(u) - \frac{1}{2k} \|u\|^2$ is strongly concave.

Let $\{k_s\}_{s=0}^\infty$ be a positive monotone increasing sequence and $\lim_{s \to \infty} k_s = \infty$. We call it a regularization sequence. The correspondent sequence $\{\lambda_s\}_{s=0}^\infty$:

$$(35) \qquad\qquad \lambda_s = \operatorname{argmax}\{d(u) - \frac{1}{2k_s} \|u\|^2 : u \in \mathbb{R}^m\}$$

is unique due to the strong concavity of $D(u, k)$ in $u$.

**Theorem 3.** *If* $L^* = \operatorname{Argmax}\{d(\lambda) | \lambda \in \mathbb{R}^m\}$ *is bounded and* $f$, $c_i \in C^1$, $i = 1, ..., m$, *then for any regularization sequence* $\{k_s\}_{s=0}^\infty$ *the following statements hold*

1) $\|\lambda_{s+1}\| > \|\lambda_s\|$;
2) $d(\lambda_{s+1}) > d(\lambda_s)$;
3) $\lim_{s \to \infty} \lambda_s = \lambda^* = \operatorname{argmin}_{\lambda \in L^*} \|\lambda\|$.

*Proof.* It follows from (35) and strong concavity of $D(u, k)$ in $u \in \mathbb{R}^m$ that

$$d(\lambda_s) - (2k_s)^{-1} \|\lambda_s\|^2 > d(\lambda_{s+1}) - (2k_s)^{-1} \|\lambda_{s+1}\|^2$$

and

$$(36) \qquad d(\lambda_{s+1}) - (2k_{s+1})^{-1} \|\lambda_{s+1}\|^2 > d(\lambda_s) - (2k_{s+1})^{-1} \|\lambda_s\|^2.$$

By adding the inequalities we obtain

$$(37) \qquad\qquad 0.5(k_s^{-1} - k_{s+1}^{-1})[\|\lambda_{s+1}\|^2 - \|\lambda_s\|^2] > 0.$$

Keeping in mind $k_{s+1} > k_s$ from (37) we obtain 1).

From (36) we have

$$(38) \qquad d(\lambda_{s+1}) - d(\lambda_s) > (2k_{s+1})^{-1}[\|\lambda_{s+1}\|^2 - \|\lambda_s\|^2] > 0,$$

therefore from 1) follows 2).

Due to concavity $d$ from boundedness of $L^*$ follows boundedness of any level set $\Lambda(\lambda_0) = \{\lambda \in \mathbb{R}^m : d(\lambda) \geq d(\lambda_0)\}$ (see Theorem 24 [22]). From 2) follows $\{\lambda_s\}_{s=0}^{\infty} \subset \Lambda(\lambda_0)$, therefore for any converging subsequence $\{\lambda_{s_i}\} \subset \{\lambda_s\}_{s=0}^{\infty}$: $\lim_{s_i \to \infty} \lambda_{s_i} = \hat{\lambda}$ we have

$$(39) \qquad d(\lambda_{s_i}) - (2k_{s_i})^{-1}\|\lambda_{s_i}\|^2 > d(\lambda^*) - (2k_{s_i})^{-1}\|\lambda^*\|^2.$$

Taking the limit in (39) when $k_{s_i} \to \infty$ we obtain $d(\hat{\lambda}) \geq d(\lambda^*)$, therefore $\hat{\lambda} = \lambda^* \in L$. In view of 2) we have $\lim_{s \to \infty} d(\lambda_s) = d(\lambda^*)$.

It follows from 1) that $\lim_{s \to \infty} \|\lambda_s\| = \|\lambda^*\|$. Also from

$$d(\lambda_s) - (2k_s)^{-1}\|\lambda_s\|^2 > d(\lambda^*) - (2k_s)^{-1}\|\lambda^*\|^2$$

follows

$$\|\lambda^*\|^2 - \|\lambda_s\|^2 > 2k_s(d(\lambda^*) - d(\lambda_s)) \geq 0, \quad \forall \lambda^* \in L^*,$$

therefore $\lim_{s \to \infty} \|\lambda_s\| = \min_{\lambda \in L^*} \|\lambda\|$.

Convergence of the regularization method (34) is due to unbounded increase of the penalty parameter $k > 0$, therefore one can hardly expect solving the problem (23) with high accuracy.

## 3.2. Augmented Lagrangian and Quadratic Proximal Point Method.
In this subsection we consider Augmented Lagrangian method (see [28], [52]), which allows eliminate difficulties associated with unbounded increase of the penalty parameter.

The problem (16) is equivalent to the following problem

$$(40) \qquad f(x) + k^{-1} \sum_{i=1}^{m} \pi(kc_i(x)) \to \min$$

$$(41) \qquad \text{s.t. } c_i(x) = 0, \ \ i = 1, ..., m.$$

The correspondent classical Lagrangian $\mathcal{L} : \mathbb{R}^n \times \mathbb{R}^m \times \mathbb{R}_{++} \to \mathbb{R}$ for the equivalent problem (40)-(41) is given by

$$\mathcal{L}(x, \lambda, k) = f(x) - \sum_{i=1}^{m} \lambda_i c_i(x) + k^{-1} \sum_{i=1}^{m} \pi(kc_i(x)) =$$

$$f(x) - \sum_{i=1}^{m} \lambda_i c_i(x) + \frac{k}{2} \sum_{i=1}^{m} c_i^2(x).$$

$\mathcal{L}$ is called Augmented Lagrangian (AL) for the original problem (16).

We assume that for a given $(\lambda, k) \in \mathbb{R}^m \times \mathbb{R}_{++}^1$ the unconstrained minimizer $\hat{x}$ exists, that is

$$(42) \qquad \hat{x} = \hat{x}(\lambda, k) : \nabla_x \mathcal{L}(\hat{x}, \lambda, k) = \nabla f(\hat{x}) - \sum_{i=1}^{m} (\lambda_i - \pi'(kc_i(\hat{x}))) \nabla c_i(\hat{x}) = 0.$$

Let

$$(43) \qquad \hat{\lambda}_i = \hat{\lambda}_i(\lambda, k) = \lambda_i - \pi'(kc_i(\hat{x})), i = 1, ..., m.$$

Then from (42) follows $\nabla_x L(\hat{x}, \hat{\lambda}) = 0$, which means that $\hat{x}$ satisfies the necessary condition for $\hat{x}$ to be a minimizer of $L(x, \hat{\lambda})$. If $L(\hat{x}, \hat{\lambda}) = \min_{x \in \mathbb{R}^n} L(x, \hat{\lambda})$ then $d(\hat{\lambda}) = L(\hat{x}, \hat{\lambda})$ and

$$(44) \qquad\qquad\qquad c(\hat{x}) \in \partial d(\hat{\lambda}).$$

From (43) follows

$$(45) \qquad\qquad\qquad c(\hat{x}) = \frac{1}{k}\pi^{'-1}(\hat{\lambda} - \lambda).$$

Using LEID and (45) we obtain

$$0 \in \partial d(\hat{\lambda}) - k^{-1}\sum_{i=1}^{m}\pi^{*'}(\hat{\lambda}_i - \lambda)e_i,$$

which is the optimality condition for $\hat{\lambda}$ to be the maximizer in the following unconstrained maximization problem

$$(46) \quad d(\hat{\lambda}) - k^{-1}\sum_{i=1}^{m}\pi^*(\hat{\lambda}_i - \lambda_i) = \max\{d(u) - k^{-1}\sum_{i=1}^{m}\pi^*(u_i - \lambda_i) : u \in \mathbb{R}^n\}.$$

In view of $\pi^*(s) = \frac{1}{2}s^2$ we can rewrite (46) as follows

$$(47) \qquad\qquad \hat{\lambda} = \operatorname{argmax}\{d(u) - \frac{1}{2k}\|u - \lambda\|^2 : u \in \mathbb{R}^n\}$$

Thus the multipliers method (42)-(43) is equivalent to the quadratic proximal point (prox) method (47) for the dual problem (23) (see [27],[35], [38], [56]-[58] and references therein)

If $\hat{x}$ is a unique solution to the system $\nabla_x L(x, \hat{\lambda}) = 0$, then $\nabla d(\hat{\lambda}) = c(\hat{x})$ and from (45) follows

$$\hat{\lambda} = \lambda + k\nabla d(\hat{\lambda}),$$

which is an implicit Euler method for solving the following system of ordinary differential equations

$$(48) \qquad\qquad\qquad \frac{d\lambda}{dt} = k\nabla d(\lambda), \ \lambda(0) = \lambda_0.$$

Let us consider the prox-function $p : \mathbb{R}^m \to \mathbb{R}$ defined as follows

$$p(\lambda) = d(u(\lambda)) - \frac{1}{2k}\|u(\lambda) - \lambda\|^2 = D(u(\lambda), \lambda) =$$

$$\max\{d(u) - \frac{1}{2k}\|u - \lambda\|^2 : u \in \mathbb{R}^n\}.$$

The function $D(u, \lambda)$ is strongly concave in $u \in \mathbb{R}^m$, therefore $u(\lambda) = \operatorname{argmax}\{D(u, \lambda) : u \in \mathbb{R}^n\}$ is unique. The prox-function $p$ is concave and differentiable. For its gradient we have

$$\nabla p(\lambda) = \nabla_u D(u(\lambda), \lambda) \cdot \nabla_\lambda u(\lambda) + \nabla_\lambda D(u(\lambda), \lambda),$$

where $\nabla_\lambda u(\lambda)$ is the Jacobian of $u(\lambda) = (u_1(\lambda), ..., u_m(\lambda))^T$. Keeping in mind $\nabla_u D(u(\lambda), \lambda) = 0$ we obtain

$$\nabla p(\lambda) = \nabla_\lambda D(u, \lambda) = \frac{1}{k}(u(\lambda) - \lambda) = \frac{1}{k}(\hat{\lambda} - \lambda)$$

or

$$(49) \qquad\qquad\qquad \hat{\lambda} = \lambda + k\nabla p(\lambda).$$

In other words, the prox-method (47) is an explicit Euler method for the following system

$$\frac{d\lambda}{dt} = k\nabla p(\lambda), \ \lambda(0) = \lambda_0.$$

By reiterating (49) we obtain the dual sequence $\{\lambda_s\}_{s=0}^{\infty}$:

(50) $$\lambda_{s+1} = \lambda_s + k\nabla p(\lambda_s),$$

generated by the gradient method for maximization the prox function $p$. The gradient $\nabla p$ satisfies Lipschitz condition with constant $L = k^{-1}$. Therefore we have the following bound $\Delta p(\lambda_s) = p(\lambda^*) - p(\lambda_s) \leq O(sk)^{-1}$ (see, for example, [45]).

We saw that the dual aspects of the penalty and the multipliers methods are critical for understanding their convergence properties and LEID is the main instrument for obtaining the duality results.

It is even more so for constrained optimization problems with inequality constraints.

## 4. SUMT as Interior Regularization Methods for the Dual Problem

The sequential unconstrained minimization technique (SUMT) (see [22]) goes back to the 50s, when R.Frisch introduced log-barrier function to replace a convex optimization with inequality constraints by a sequence of unconstrained convex minimization problems.

Let $f$ and all-$c_i$, $i = 1, ..., m$ be convex and smooth. We consider the following convex optimization problem

(51) $$\min f(x)$$
$$\text{s. t. } x \in \Omega,$$

where $\Omega = \{x \in \mathbb{R}^n : c_i(x) \geq 0, \ i = 1, ..., m\}$.

From this point on we assume

A. The solution set $X^* = \text{Argmin}\{f(x) : x \in \Omega\}$ is not empty and bounded.

B. Slater condition holds, i.e. there exists $x_0 \in \Omega$: $c_i(x_0) > 0$, $i = 1, ..., m$.

By adding one constraint $c_0(x) = M - f(x) \geq 0$ with $M$ large enough to the original set of constraints $c_i(x) \geq 0$, $i = 1, ..., m$ we obtain a new feasible set, which due to the assumption A convexity $f$ and concavity $c_i$, $i = 1, ..., m$ is bounded (see Theorem 24 [22]) and the extra constraint $c_0(x) \geq 0$ for large $M$ does not effect $X^*$.

So we assume from now on that $\Omega$ is bounded. It follows from KKT's Theorem that under Slater condition the existence of the primal solution

$$f(x^*) = \min\{f(x)|x \in \Omega\}$$

leads to the existence of $\lambda^* \in \mathbb{R}_+^m$ that for $\forall x \in \mathbb{R}^n$ and $\lambda \in \mathbb{R}_+^m$ we have

(52) $$L(x^*, \lambda) \leq L(x^*, \lambda^*) \leq L(x, \lambda^*)$$

and $\lambda^*$ is the solution of the dual problem

(53) $$d(\lambda^*) = \max\{d(\lambda)|\lambda \in \mathbb{R}_+^m\}.$$

Also from B follows boundedness of the dual optimal set

$$L^* = \text{Argmax}\{d(\lambda) : \lambda \in \mathbb{R}_+^m\}.$$

From concavity $d$ and boundedness $L^*$ follows boundedness of the dual level set $\Lambda(\bar{\lambda}) = \{\lambda \in \mathbb{R}^m_+ : d(\lambda) \geq d(\bar{\lambda})\}$ for any given $\bar{\lambda} \in \mathbb{R}^m_+$: $d(\bar{\lambda}) < d(\lambda^*)$.

4.1. **Logarithmic Barrier.** To replace the constrained optimization problem (51) by a sequence of unconstrained minimization problems R. Frisch in 1955 introduced (see [23]) the log-barrier penalty function $P : \mathbb{R}^n \times \mathbb{R}_{++} \to \mathbb{R}$ defined as follows

$$P(x, k) = f(x) - k^{-1} \sum_{i=1}^{m} \pi(kc_i(x)),$$

where $\pi(t) = \ln t$, $(\pi(t) = -\infty$ for $t \leq 0)$ and $k > 0$. Due to convexity $f$ and concavity $c_i$ $i = 1, ..., m$ the function $P$ is convex in $x$. Due to Slater condition, convexity $f$, concavity $c_i$ and boundedness $\Omega$ the recession cone of $\Omega$ is empty that is for any $x \in \Omega$, $k > 0$ and $0 \neq d \in \mathbb{R}^n$ we have

$$\tag{54} \lim_{t \to \infty} P(x + td, k) = \infty.$$

Therefore for any $k > 0$ there exists

$$\tag{55} x(k) : \nabla_x P(x(k), k) = 0.$$

**Theorem 4.** *If $A$ and $B$ hold and $f$, $c_i \in C^1$, $i = 1, ..., m$, then interior log-barrier method* (55) *is equivalent to the interior regularization method*

$$\tag{56} \lambda(k) = \operatorname{argmax}\{d(u) + k^{-1} \sum_{i=1}^{m} \ln u_i : u \in \mathbb{R}^m_+\}$$

*and the following error bound holds*

$$\tag{57} \max\{\Delta f(x(k)) = f(x(k)) - f(x^*), \Delta d(\lambda(k)) = d(\lambda^*) - d(\lambda(k))\} = mk^{-1}.$$

*Proof.* From (54) follows existence $x(k) : P(x(k), k) = \min\{P(x, k) : x \in \mathbb{R}^n\}$ for any $k > 0$.

Therefore

$$\tag{58} \nabla_x P(x(k), k) = \nabla f(x(k)) - \sum_{i=1}^{m} \pi^{'}(k_i(x(k)) \nabla c_i(x(k)) = 0.$$

Let

$$\tag{59} \lambda_i(k) = \pi^{'}(kc_i(x(k)) = (kc_i(x(k)))^{-1}, \ i = 1, .., m.$$

Then from (58) and (59) follows $\nabla_x P(x(k), k) = \nabla_x L(x(k), \lambda(k)) = 0$, therefore $d(\lambda(k)) = L(x(k), \lambda(k))$. From $\pi^{''}(t) = -t^2 < 0$ follows existence of $\pi^{'-1}$ and from (59) we have $kc(x(k)) = \pi^{'-1}(\lambda_i(k))$. Using LEID we obtain

$$\tag{60} c_i(x(k)) = k^{-1} \pi^{*'}(\lambda_i(k)),$$

where $\pi^*(s) = \inf_{t>0}\{st - \ln t\} = 1 + \ln s$. The subdifferential $\partial d(\lambda(k))$ contains $-c(x(k))$, that is

$$\tag{61} 0 \in \partial d(\lambda(k)) + c(x(k)).$$

From (60) and (61) follows

$$\tag{62} 0 \in \partial d(\lambda(k)) + k^{-1} \sum_{i=1}^{m} \pi^{*'}(\lambda_i(k))e_i.$$

The last inclusion is the optimality criteria for $\lambda(k)$ to be the maximizer in (56).

The maximizer $\lambda(k)$ is unique due to the strict concavity of the objective function in (56).

Thus, SUMT with log-barrier function $P(x, k)$ is equivalent to the interior regularization method (56).

For primal interior trajectory $\{x(k)\}_{k=k_0>0}^{\infty}$ and dual interior trajectory $\{\lambda(k)\}_{k=k_0>0}^{\infty}$ we have

$$f(x(k)) \geq f(x^*) = d(\lambda^*) \geq d(\lambda(k)) = L(x(k), \lambda(k)) = f(x(k)) - (c(x(k)), \lambda(k)).$$

From (59) follows $\lambda_i(k) c_i(x(k)) = k^{-1}$, $i = 1, ..., m$, hence for the primal-dual gap we obtain

$$f(x(k)) - d(\lambda(k)) = (c(x(k)), \lambda(k)) = m k^{-1}.$$

Therefore for the primal and the dual error bounds we obtain (57).          □

The main idea of the interior point methods (IPMs) is to stay "close" to the primal $\{x(k)\}_{k=0}^{\infty}$ or to the primal-dual $\{x(k), \lambda(k)\}_{k=0}^{\infty}$ trajectory and increase $k > 0$ at each step by a factor $(1 - \frac{\alpha}{\sqrt{n}})^{-1}$, where $\alpha > 0$ is independent of $n$. In case of LP at each step we solve a system of linear equations, which requires $O(n^{2.5})$ operations. Therefore accuracy $\varepsilon > 0$ IPM are able to achieve in $O(n^3 \ln \varepsilon^{-1})$ operations.

In case of log-barrier transformation the situation is symmetric, that is both the primal interior penalty method (55) and the dual interior regularization method (56) are using the same log-barrier function.

It is not the case for other constraints transformations used in SUMT.

## 4.2. Hyperbolic Barrier.
The hyperbolic barrier

$$\pi(t) = \begin{cases} -t^{-1}, t > 0 \\ -\infty, t \leq 0, \end{cases}$$

has been introduced by C. Carroll in the 60s, (see [12]). It leads to the following hyperbolic penalty function

$$P(x, k) = f(x) - k^{-1} \sum_{i=1}^{m} \pi(k c_i(x)) = f(x) + k^{-1} \sum_{i=1}^{m} (k c_i(x))^{-1},$$

which is convex in $x \in \mathbb{R}^n$ for any $k > 0$. For the primal minimizer we obtain

$$(63) \qquad x(k) : \nabla_x P(x(k), k) = \nabla f(x(k)) - \sum_{i=1}^{m} \pi^{'}(k c_i(x(k))) \nabla c_i(x(k)) = 0.$$

For the vector of Lagrange multipliers we have

$$(64) \qquad \lambda(k) = (\lambda_i(k) = \pi^{'}(k c_i(x(k))) = (k c_i(x(k)))^{-2}, \ i = 1, ..., m).$$

We will show later that vectors $\lambda(k)$, $k \geq 1$ are bounded. Let $L = \max_{i,k} \lambda_i(k)$.

**Theorem 5.** *If A and B hold and $f$, $c_i \in C^1$, $i = 1, .., m$, then hyperbolic barrier method (63) is equivalent to the parabolic regularization method*

$$(65) \qquad d(\lambda(k)) + 2k^{-1} \sum_{i=1}^{m} \sqrt{\lambda_i(k)} = \max\{d(u) + 2k^{-1} \sum_{i=1}^{m} \sqrt{u_i} : u \in \mathbb{R}_+^m\}$$

*and the following bounds holds*

$$\max\{\Delta f(x(k)) = f(x(k)) - f(x^*),$$

(66)                    $\Delta d(\lambda(k)) = d(\lambda^*) - d(\lambda(k))\} \leq m\sqrt{L}k^{-1}.$

*Proof.* From (63) and (64) follows

$$\nabla_x P(x(k), k) = \nabla_x L(x(k), \lambda(k)) = 0,$$

therefore $d(\lambda(k)) = L(x(k), \lambda(k)).$

From $\pi''(t) = -2t^{-3} < 0, \forall t > 0$ follows existence of $\pi'^{-1}.$

Using LEID from (64) we obtain

$$c_i(x(k)) = k^{-1}\pi'^{-1}(\lambda_i(k)) = k^{-1}\pi^{*'}(\lambda_i(k)), \ i = 1, ..., m,$$

where $\pi^*(s) = \inf_t\{st - \pi(t)\} = 2\sqrt{s}.$

The subgradient $-c(x(k)) \in \partial d(\lambda(k))$ that is

(67)        $0 \in \partial d(\lambda(k)) + c(x(k)) = \partial d(\lambda(k)) + k^{-1} \sum_{i=1}^{m} \pi^{*'}(\lambda_i(k))e_i.$

The last inclusion is the optimality condition for the interior regularization method (65) for the dual problem.

Thus, the hyperbolic barrier method (63) is equivalent to the parabolic regularization method (65) and $D(u, k) = d(u) + 2k^{-1} \sum_{i=1}^{m} \sqrt{u_i}$ is strictly concave.

Using considerations similar to those in Theorem 3 and keeping in mind strict concavity of $D(u, k)$ in $u$ from (65) we obtain

$$\sum_{i=1}^{m} \sqrt{\lambda_i(1)} > ... \sum_{i=1}^{m} \sqrt{\lambda_i(k)} > \sum_{k=1}^{m} \sqrt{\lambda_i(k+1)} > ...$$

Therefore the sequence $\{\lambda(k)\}_{k=1}^{\infty}$ is bounded, so there exists $L = \max_{i,k} \lambda_i(k) > 0.$ From (64) for any $k \geq 1$ and $i = 1, ..., m$ we have

$$\lambda_i(k)c_i^2(x(k)) = k^{-2}$$

or

$$(\lambda_i(k)c_i(x(k)))^2 = k^{-2}\lambda_i(k) \leq k^{-2}L.$$

Therefore

$$(\lambda(k), c(x(k))) \leq m\sqrt{L}k^{-1}.$$

For the primal interior sequence $\{x(k)\}_{k=1}^{\infty}$ and dual interior sequence $\{\lambda(k)\}_{k=1}^{\infty}$ we have

$$f(x(k)) \geq f(x^*) = d(\lambda^*) \geq L(x(k), \lambda(k)) = d(\lambda(k)),$$

therefore

$$f(x(k)) - d(\lambda(k)) = (c(x(k)), \lambda(k))) \leq m\sqrt{L}k^{-1},$$

which leads to (66).                                                                     □

In spite of similarity bounds (57) and (65) are fundamentally different because $L$ can be very large for problems where Slater condition is "barely" satisfied, that is the primal feasible set is not "well defined".

This is one of the reasons why log-barrier function is so important.

4.3. **Exponential Penalty.** Exponential penalty $\pi(t) = -e^{-t}$ has been used by Motzkin in 1952 (see [40]) to transform a systems of linear inequalities into an unconstrained convex optimization problem in order to use unconstrained minimization technique for solving linear inequalities.

The exponential transformation $\pi(t) = -e^{-t}$ leads to the exponential penalty function

$$P(x,k) = f(x) - k^{-1} \sum_{i=1}^{m} \pi(kc_i(x)) = f(x) + k^{-1} \sum_{i=1}^{m} e^{-kc_i(x)},$$

which is for any $k > 0$ convex in $x \in \mathbb{R}^n$.

For the primal minimizer we have

$$(68) \qquad x(k): \nabla_x P(x(k),k) = \nabla f(x(k)) - \sum_{i=1}^{m} e^{-kc_i(x(k))} \nabla c_i(x(k)) = 0.$$

Let us introduce the Lagrange multipliers vector

$$(69) \qquad \lambda(k) = (\lambda_i(k) = \pi'(c_i(x(k)) = e^{-kc_i(x(k))}, \ i = 1,...,m)$$

From (68) and (69) we have

$$\nabla_x P(x(k),k) = \nabla_x L(x(k), \lambda(k)) = 0.$$

Therefore from convexity $L(x, \lambda(k))$ in $x \in \mathbb{R}^n$ follows $d(\lambda(k)) = \min\{L(x, \lambda(k)) | x \in \mathbb{R}^n\} = L(x(k), \lambda(k))$ and $-c(x(k)) \in \partial d(\lambda(k))$, therefore

$$(70) \qquad 0 \in c(x(k)) + \partial d(\lambda(k)).$$

From $\pi''(t) = -e^{-t} \neq 0$ follows the existence $\pi'^{-1}$, therefore using LEID from (69) we obtain

$$c_i(x(k)) = k^{-1}\pi'^{-1}(\lambda_i(k)) = k^{-1}\pi^{*'}(\lambda_i(k)), \ i = 1,...,m.$$

Inclusion (70) we can rewrite as follows

$$\partial d(\lambda(k)) + k^{-1} \sum \pi^{*'}(\lambda(k))e_i = 0.$$

Keeping in mind $\pi^*(s) = \inf_t\{st - \pi(t)\} = \inf\{st + e^{-t}\} = -s \ln s + s$ from the last inclusion we obtain

$$d(\lambda(k)) - k^{-1} \sum_{i=1}^{m} \lambda_i(k)(\ln(\lambda_i(k) - 1)) =$$

$$(71) \qquad \max\{d(u) - k^{-1} \sum_{i=1}^{m} u_i(\ln u_i - 1) : u \in \mathbb{R}_+^m\}.$$

It means that the exponential penalty method (68) is equivalent to the interior regularization method (71) with strictly concave regularization function $r(u) = -\sum_{i=1}^{m} u_i(\ln u_i - 1)$.

The convergence of the dual sequence $\{\lambda(k)\}_{k=0}^{\infty}$ can be proven using arguments similar to those used in Theorem 3.

We conclude the section by considering smoothing technique for convex optimization.

4.4. **Log-Sigmoid (LS) Method.** It follows from Karush-Kuhn-Tucker's Theorem that under Slater condition for $x^*$ to be a solution of (51) it is necessary and sufficient existence $\lambda^* \in \mathbb{R}^m$, that the pair $(x^*; \lambda^*)$ is the saddle point of the Lagrangian, that is (52) hold.

From the right inequality of (52) and complementarity condition we obtain

$$f(x^*) \leq f(x) - \sum_{i=1}^m \lambda_i^* \min\{c_i(x), 0\} \leq$$

$$f(x) - \max_{1 \leq i \leq m} \lambda_i^* \sum_{i=1}^m \min\{c_i(x), 0\}$$

for any $x \in \mathbb{R}^n$. Therefore for any $r > \max_{1 \leq i \leq m} \lambda_i^*$ we have

$$(72) \qquad f(x^*) \leq f(x) - r \sum_{i=1}^m \min\{c_i(x), 0\}, \forall x \in \mathbb{R}^n.$$

The function

$$Q(x, r) = f(x) - r \sum_{i=1}^m \min\{c_i(x), 0\}$$

is called exact penalty function.

Due to concavity $c_i$, $i = 1, ..., m$ functions $q_i(x) = \min\{c_i(x), 0\}$ are concave. From convexity $f$ and concavity $q_i$, $i = 1, ..., m$ follows convexity $Q(x, r)$ in $x \in \mathbb{R}^n$. From (72) follows that solving (51) is equivalent to solving the following unconstrained minimization problem

$$(73) \qquad f(x^*) = Q(x^*, r) = \min\{Q(x, r) : x \in \mathbb{R}^n\}.$$

The function $Q(x, r)$ is non-smooth at $x^*$. The smoothing techniques replace $Q$ by a sequence of smooth functions, which approximate $Q(x, r)$. (see [14], [47], [48] and references therein)

Log-sigmoid (LS) function $\pi : \mathbb{R} \to \mathbb{R}$ is defined by

$$\pi(t) = \ln S(t, 1) = \ln(1 + e^{-t})^{-1},$$

is one of such functions. We collect the log-sigmoid properties in the following assertion

*Assertion* 1. The following statements are holds

1. $\pi(t) = t - \ln(1 + e^t) < 0$, $\pi(0) = -\ln 2$
2. $\pi^{'}(t) = (1 + e^t)^{-1} > 0$, $\pi^{'}(0) = 2^{-1}$
3. $\pi^{''}(t) = -e^t(1 + e^t)^{-2} < 0$, $\pi^{''}(0) = -2^{-2}$.

The smooth penalty method employs the scaled LS function

$$(74) \qquad k^{-1}\pi(kt) = t - k^{-1}\ln(1 + e^{kt}),$$

which is a smooth approximation of $q(t) = \min\{t, 0\}$.

In particular, from (74) follows

$$(75) \qquad 0 < q(t) - k^{-1}\pi(kt) < k^{-1}\ln 2.$$

It means that by increasing $k > 0$ the approximation can be made as accurate as one wants.

The smooth penalty function $P : \mathbb{R}^n \times \mathbb{R}_{++} \to \mathbb{R}$ defined by

$$(76) \qquad P(x, k) = f(x) - k^{-1} \sum_{i=1}^{m} \pi(kc_i(x))$$

is the main instrument in the smoothing technique.

From Assertion 1 follows that $P$ is as smooth as $f$ and $c_i$, $i = 1, ..., m$.

The LS method at each step finds

$$(77) \qquad x(k) : P(x(k), k) = \min\{P(x, k) : x \in \mathbb{R}^n\}$$

and increases $k > 0$ if the accuracy obtained is not satisfactory.

Without loss of generality we assume that $f$ is bounded from below. Such assumption does not restrict the generality, because the original objective function $f$ can be replaced by an equivalent $f(x) := \ln(1 + e^{f(x)}) \geq 0$.

Boundedness of $\Omega$ together with Slater condition, convexity $f$ and concavity $c_i$, $i = 1, ..., m$ make the recession cone of $\Omega$ empty, that is (54) holds for $P(x, k)$ given by (76), any $k > 0$, $d \in \mathbb{R}^n$ and any $x \in \Omega$.

Therefore minimizer $x(k)$ in (77) exists for any $k > 0$ that is

$$\nabla_x P(x(k), k) = \nabla f(x(k)) - \sum_{i=1}^{m} \pi^{'}(kc_i(x(k)))\nabla c_i(x(k)) =$$

$$= \nabla f(x(k)) - \sum_{i=1}^{m} (1 + e^{kc_i(x(k))})^{-1}\nabla c_i(x(k)) = 0.$$

Let

$$(78) \qquad \lambda_i(k) = (1 + e^{kc_i(x(k))})^{-1}, \ i = 1, ..., m,$$

then

$$\nabla_x P(x(k); k) = \nabla f(x(k)) - \sum_{i=1}^{m} \lambda_i(k)\nabla c_i(x(k)) = 0.$$

From (78) follows $\lambda_i(k) \leq 1$ for any $k > 0$. Therefore, generally speaking, one can't expect finding a good approximation for optimal Lagrange multipliers, no matter how large the penalty parameter $k > 0$ is.

If the dual sequence $\{\lambda(k)\}_{k=k_0}^{\infty}$ does not converges to $\lambda^* \in L^*$, then in view of the last equation one can't expect convergence of the primal sequence $\{x(k)\}_{k=k_0}^{\infty}$ to $x^* \in X^*$.

To guarantee convergence of the LS method we have to modify $P(x, k)$. Let $0 < \alpha < 0.5$ and

$$(79) \qquad P(x, k) := P_\alpha(x, k) = f(x) - k^{-1+\alpha} \sum_{i=1}^{m} \pi(kc_i(x)).$$

It is easy to see that the modification does not effect the existence of $x(k)$. Therefore for any $k > 0$ there exists

$$(80) \qquad x(k) : \nabla_x P(x(k), k) = \nabla f(x(k)) - k^{\alpha} \sum \pi^{'}(kc(x(k)))\nabla c_i(x(k)) = 0.$$

**Theorem 6.** *If A and B hold and $f$, $c_i \in C^1$, $i = 1, ..., m$, then the LS method (80) is equivalent to an interior regularization method*

$$d(\lambda(k)) + k^{-1} \sum_{i=1}^{m} \pi^*(k^{-\alpha}\lambda_i(k)) =$$

$$\max\{d(u) + k^{-1} \sum_{i=1}^{m} \pi^*(k^{-\alpha}u_i) : 0 \le u_i \le k^\alpha, \ i = 1, ..., m\}.$$

*Proof.* Let

$$(81) \qquad \lambda_i(k) = k^\alpha \pi^{'}(kc_i(x(k))) = k^\alpha(1 + e^{kc_i(x(k))})^{-1}, \ i = 1, ..., m.$$

From (80) and (81) follows

$$(82) \qquad \begin{aligned} \nabla_x P(x(k), k) &= \nabla f(x(k)) - \sum_{i=1}^{m} \lambda_i(k) \nabla c_i(x(k)) = \\ &\nabla_x L(x(k), \lambda(k)) = 0. \end{aligned}$$

From (81) we have

$$(83) \qquad \pi^{'}(kc_i(x(k))) = k^{-\alpha}\lambda_i(k).$$

Due to $\pi^{''}(t) < 0$ there exists $\pi^{'-1}$, therefore

$$c_i(x(k)) = k^{-1}\pi^{'-1}(k^{-\alpha}\lambda_i(k)).$$

Using LEID we obtain

$$(84) \qquad c_i(x(k)) = k^{-1}\pi^{*'}(k^{-\alpha}\lambda_i(k)),$$

where

$$\pi^*(s) = \inf_t \{st - \pi(t)\} = -[(1-s)\ln(1-s) + s\ln s]$$

is Fermi-Dirac (FD) entropy function (see, for example, [54]).

From (82) follows $d(\lambda(k)) = L(x(k), \lambda(k))$ , also the subdifferential $\partial d(\lambda(k))$ contains $-c(x(k))$, that is

$$(85) \qquad 0 \in c(x(k)) + \partial d(\lambda(k)).$$

Combining (84) and (85) we obtain

$$(86) \qquad 0 \in \partial d(\lambda(k)) + k^{-1} \sum_{i=1}^{m} \pi^{*'}(k^{-\alpha}\lambda_i(k))e_i.$$

The inclusion (86) is the optimality criteria for the following problem

$$d(\lambda(k)) + k^{-1} \sum_{i=1}^{m} \pi^*(k^{-\alpha}\lambda_i(k)) =$$

$$(87) \qquad \max\{d(u) + k^{-1}r(u) : 0 \le u_i \le k^\alpha, \ i = 1, .., m\},$$

where $r(u) = \sum_{i=1}^{m} \pi^*(k^{-\alpha}u_i)$.

In other words the LS method (80)-(81) is equivalent to the interior regularization method (87) with FD entropy function used for dual regularization. The FD function is strongly concave inside the cube $\{u \in \mathbb{R}^m : 0 \le u_i \le k^\alpha, i = 1, ..., m\}$.

It follows from (87) that for any regularization sequence $\{k_s\}_{s=0}^\infty$ the Lagrange multipliers $0 < \lambda_i(k_s) < k_s^\alpha$, $i = 1, ..., m$ can be any positive number, which underlines the importance of modification (79).

**Theorem 7.** *Under conditions of Theorem 6 for any regularization sequence $\{k_s\}_{s=0}^{\infty}$, the primal sequence*

$$(88) \qquad \{x_s\}_{s=0}^{\infty} : \nabla_x P(x_s, k_s) = \nabla f(x_s) - \sum_{i=1}^{m} \lambda_{i,s} \nabla c_i(x_s) = 0$$

*and the dual sequence*

$$\{\lambda_s\}_{s=0}^{\infty} : d(\lambda_s) + k_s^{-1} r(\lambda_s) =$$

$$(89) \qquad \max\{d(u) + k_s^{-1} r(u) : 0 \le u_i \le k^{\alpha}, \ i = 1, ..., m\}$$

*the following statements hold*

1) *a)* $d(\lambda_{s+1}) > d(\lambda_s)$*; b)* $r(\lambda_{s+1}) < r(\lambda_s)$*;*
2) $\lim_{s \to \infty} d(\lambda_s) = d(\lambda^*)$ *and* $\lambda^* = \arg\min\{r(\lambda) : \lambda \in L^*\}$*;*
3) *the primal-dual sequence* $\{x_s, \lambda_s\}_{s=0}^{\infty}$ *is bounded and any limit point is the primal-dual solution.*

*Proof.*     1) From (89) and strong concavity $r(u)$ follows

$$(90) \qquad d(\lambda_{s+1}) + k_{s+1}^{-1} r(\lambda_{s+1}) > d(\lambda_s) + k_{s+1}^{-1} r(\lambda_s)$$

and

$$(91) \qquad d(\lambda_s) + k_s^{-1} r(\lambda_s) > d(\lambda_{s+1}) + k_s^{-1} r(\lambda_{s+1}).$$

Therefore

$$(k_{s+1}^{-1} - k_s^{-1})(r(\lambda_{s+1}) - r(\lambda_s)) > 0.$$

From $k_{s+1} > k_s$ and last inequality follows $r(\lambda_{s+1}) < r(\lambda_s)$, therefore from (90) follows

$$(92) \qquad d(\lambda_{s+1}) > d(\lambda_s) + k_{s+1}^{-1}(r(\lambda_s) - r(\lambda_{s+1})) > d(\lambda_s).$$

2) The monotone increasing sequence $\{d(\lambda_s)\}_{s=0}^{\infty}$ is bounded from above by $f(x^*)$. Therefore there is $\lim_{s \to \infty} d(\lambda_s) = \bar{d} \le f(x^*) = d(\lambda^*)$.

From (89) follows

$$(93) \qquad d(\lambda_s) + k_s^{-1} r(\lambda_s) \ge d(\lambda^*) + k_s^{-1} r(\lambda^*).$$

From (92) follows $\{\lambda_s\}_{s=0}^{\infty} \subset \Lambda(\lambda_0) = \{\lambda \in \mathbb{R}_+^m : d(\lambda) \ge d(\lambda_0)\}$. The set $\Lambda(\lambda_0)$ is bounded due to the boundedness of $L^*$ and concavity $d$. Therefore there exists $\{\lambda_{s_i}\}_{i=1}^{\infty} \subset \{\lambda_s\}_{s=0}^{\infty}$ that $\lim_{s_i \to \infty} \lambda_{s_i} = \bar{\lambda}$. By taking the limit in the correspondent subsequence in (93) we obtain $d(\bar{\lambda}) \ge d(\lambda^*)$, that is $d(\bar{\lambda}) = d(\lambda^*)$.

From $\lim_{s \to \infty} d(\lambda_{s_i}) = d(\lambda^*)$ and 1a) follows $\lim_{s \to \infty} d(\lambda_s) = d(\lambda^*)$.

From (93) follows

$$(94) \qquad d(\lambda^*) - d(\lambda_s) \le k_s^{-1}(r(\lambda^*) - r(\lambda_s)), \ \forall \lambda^* \in L^*,$$

therefore (94) is true for $\lambda^* = \arg\min\{r(\lambda) | \lambda \in L^*\}$.

3) We saw already the dual sequence $\{\lambda_s\}_{s=0}^{\infty}$ is bounded. Let us show that the primal is bounded too. For a given approximation $x_s$ let consider two sets of indices $I_+(x_s) = \{i : c_i(x_s) \ge 0\}$ and $I_-(x_s) = \{i : c_i(x_s) < 0\}$.

Then keeping in mind $f(x_s) \geq 0$ we obtain

$$
\begin{aligned}
(95) \quad P(x_s, k_s) &= f(x_s) + k_s^{-1+\alpha} \sum_{i \in I_-(x_s)} \ln(1 + e^{-k_s c_i(x_s)}) \\
&\quad + k_s^{-1+\alpha} \sum_{i \in I_+(x_s)} \ln(1 + e^{-k_s c_i(x_s)}) \\
&\geq f(x_s) - k_s^\alpha \sum_{i \in I_-(x_s)} c_i(x_s) + k_s^{-1+\alpha} \sum_{i \in I_-(x_s)} \ln(1 + e^{k_s c_i(x_s)}) \\
&\geq f(x_s) - k_s^\alpha \sum_{i \in I_-(x_s)} c_i(x_s) \geq -k_s^\alpha \sum_{i \in I_-(x_s)} c_i(x_s).
\end{aligned}
$$

On the other hand,

$$
P(x_s, k_s) \leq P(x^*, k_s) = f(x^*) - k_s^{-1+\alpha} \sum_{i=1}^m \pi(k_s c_i(x^*))
$$

$$
(96) \quad = f(x^*) + k_s^{-1+\alpha} \sum_{i=1}^m \ln(1 + e^{-k_s c_i(x^*)}) \leq f(x^*) + k_s^{-1+\alpha} m \ln 2.
$$

From (95) and (96) follows

$$
(97) \quad k_s^\alpha \sum_{i \in I_-(x_s)} |c_i(x_s)| \leq f(x^*) + k_s^{-1+\alpha} m \ln 2.
$$

Therefore for any $s \geq 1$ we have

$$
(98) \quad \max_{i \in I_-(x_s)} |c_i(x_s)| \leq k_s^{-\alpha} f(x^*) + k_s^{-1} m \ln 2.
$$

It means that the primal sequence $\{x_s\}_{s=0}^\infty$ is bounded due to the boundedness of $\Omega$. In other words, the primal-dual sequence $\{x_s, \lambda_s\}_{s=0}^\infty$ is bounded.

Let consider a converging subsequence $\{x_{s_i}, \lambda_{s_i}\}_{i=0}^\infty$: $\bar{x} = \lim_{i \to \infty} x_{s_i}$; $\bar{\lambda} = \lim_{i \to \infty} \lambda_{s_i}$. From (81) follows $\bar{\lambda}_i = 0$ for $i : c_i(\bar{x}) > 0$ and $\bar{\lambda}_i \geq 0$ for $i : c_i(\bar{x}) = 0$. From (82) follows $\nabla_x L(\bar{x}, \bar{\lambda}) = 0$, therefore $(\bar{x}, \bar{\lambda})$ is KKT's pair, that is $\bar{x} = x^*$, $\bar{\lambda} = \lambda^*$. $\qquad \square$

The equivalence primal SUMT and dual interior regularization methods not only allows to prove convergence in a unified and simple manner, but also provide important information about dual feasible solution, which can be used to improve numerical performance. One can't, however, expect finding solution with high accuracy because finding the primal minimizer for large $k > 0$ is a difficult task for the well known reasons.

The difficulties, to a large extend, one can overcome by using the Nonlinear Rescaling theory and methods (see [31], [46], [47], [50], [53], [59] and references). One can view NR as an alternative to SUMT.

## 5. Nonlinear Rescaling and Interior Prox with Entropy like Distance

The NR scheme employs smooth, strictly concave and monotone increasing functions $\psi \in \Psi$ to transform the original set of constraints into an equivalent set. The transformation is scaled by a positive scaling (penalty) parameter. The Lagrangian for the equivalent problem is our main instrument.

At each step NR finds the primal minimizer of the Lagrangian for the equivalent problem and uses the minimizer to update the Lagrange multipliers (LM). The positive scaling parameter can be fixed or updated from step to step. The fundamental difference between NR ans SUMT lies in the role of the LM vector.

In case of SUMT the LM vector is just a by product of the primal minimization. It provides valuable information about the dual vector but it does not effect the computational process. Therefore without unbound increase of the scaling parameter, which is the only tool to control the process, one can not guarantee convergence.

In the NR scheme on the top of the scaling parameter the LM vector is a critical extra tool, which controls computations.

The NR methods converges under any fixed scaling parameter, just due to the LM update (see [31], [46], [50], [53]). If one increases the scaling parameter from step to step, as SUMT does, then instead of sublinear the superlinear convergence rate can be achieved.

The interplay between Lagrangians for the original and the equivalent problems allows to show the equivalence of the primal NR method and dual proximal point method with $\varphi$-divergence entropy type distance. The kernel of the distance $\varphi = -\psi^*$, where $\psi^*$ is the LET of $\psi$. The equivalence is the key ingredient of the convergence analysis.

We consider a class $\Psi$ of smooth functions $\psi : (a, \infty) \to \mathbb{R}, \ -\infty < a < 0$ with the following properties

1) $\psi(0) = 0$; 2) $\psi^{'}(t) > 0, \ \psi(0 = 1$; 3) $\psi^{''}(t) < 0$; 4)$\lim_{t \to \infty} \psi^{'}(t) = 0$; 5) $\lim_{t \to a_+} \psi^{'}(t) = \infty$.

From 1)-3) follows

$$\Omega = \{x \in \mathbb{R}^n : c_i(x) \geq 0, \ i = 1, ..., m\} = \{x \in \mathbb{R}^n : k^{-1}\psi(kc_i(x)) \geq 0, \ i = 1, ..., m\}$$

for any $k > 0$.

Therefore (51) is equivalent to

$$
\begin{aligned}
& \min f(x) \\
& \text{s.t. } k^{-1}\psi(kc_i(x)) \geq 0, \ i = 1, ..., m.
\end{aligned}
\tag{99}
$$

The Lagrangian $\mathcal{L} : \mathbb{R}^n \times \mathbb{R}^m_+ \times \mathbb{R}_{++} \to \mathbb{R}$ for (99) is defined as follows

$$\mathcal{L}(x, \lambda, k) = f(x) - k^{-1}\sum_{i=1}^{m}\lambda_i\psi(kc_i(x)).$$

The properties of $\mathcal{L}(x, \lambda, k)$ at the KKT pair $(x^*, \lambda^*)$ we collect in the following Assertion.

*Assertion* 2. For any $k > 0$ and any KKT pair $(x^*, \lambda^*)$ the following holds

$1^\circ$ $\mathcal{L}(x^*, \lambda^*, k) = f(x^*)$

$2^\circ$ $\nabla_x\mathcal{L}(x^*, \lambda^*, k) = \nabla f(x^*) - \sum_{i=1}^{m}\lambda_i^*\nabla c_i(x^*) = \nabla_x L(x^*, \lambda^*) = 0$

$3^\circ$ $\nabla^2_{xx}\mathcal{L}(x^*, \lambda^*, k) = \nabla^2_{xx}L(x^*, \lambda^*) + k\nabla c^T(x^*)\Lambda^*\nabla c(x^*),$

where $\nabla c(x^*) = J(c(x^*))$ is the Jacobian of $c(x) = (c_1(x), ..., c_m(x))^T$ and $\Lambda^* = I \cdot \lambda^*$.

*Remark* 8. The properties $1^0 - 3^0$ show the fundamental difference between NR and SUMT. In particular, for log-barrier penalty

$$P(x, k) = f(x) - k^{-1}\sum_{i=1}^{m}\ln c_i(x)$$

neither $P$ nor its gradient or Hessian exist at the solution $x^*$. Moreover, for any given $k > 0$ we have

$$\lim_{x \to x^*} P(x, k) = \infty.$$

On the other hand, $\mathcal{L}(x, \lambda^*, k)$ is an exact smooth approximation for the non-smooth

$$F(x, x^*) = \max\{f(x) - f(x^*), -c_i(x), \ i = 1, .., m\},$$

that is, for any given $k > 0$ we have

$$\min_{x \in \mathbb{R}^n} F(x, x^*) = F(x^*, x^*) = \min_{x \in \mathbb{R}^n} (\mathcal{L}(x, \lambda^*, k) - f(x^*)) = 0.$$

5.1. **NR and Dual Prox with $\varphi$-divergence Distance.** In this subsection we consider the NR method and its dual equivalent - the prox method with $\varphi$- divergence distance for the dual problem.

Let $\psi \in \Psi$, $\lambda_0 = e = (1, ..., 1) \in \mathbb{R}_{++}^m$ and $k > 0$ are given. The NR step consists of finding the primal minimizer

$$(100) \qquad \hat{x} \equiv \hat{x}(\lambda, k) : \nabla_x \mathcal{L}(\hat{x}, \lambda, k) = 0$$

following by the Lagrange multipliers update

$$(101) \qquad \hat{\lambda} \equiv \hat{\lambda}(\lambda, k) = (\hat{\lambda}_1, ..., \hat{\lambda}_m) : \hat{\lambda}_i = \lambda_i \psi'(k c_i(\hat{x})), \ i = 1, ..., m.$$

**Theorem 9.** *If condition A and B hold and $f$, $c_i \in C^1$, $i = 1, ..., m$, then the NR method (100)-(101) is:*

1) *well defined;*
2) *equivalent to the following prox method*

$$(102) \qquad d(\hat{\lambda}) - k^{-1} D(\hat{\lambda}, \lambda) = \max\{d(u) - k^{-1} D(u, \lambda) | u \in \mathbb{R}_{++}^m\},$$

*where $D(u, \lambda) = \sum_{i=1}^m \lambda_i \varphi(u_i/\lambda_i)$ is $\varphi$-divergence distance function based on kernel $\varphi = -\psi^*$.*

*Proof.*     1) Due to the properties 1)-3) of $\psi$, convexity $f$ and concavity of all $c_i$, the Lagrangian $\mathcal{L}$ is convex in $x$. From boundedness of $\Omega$, Slater condition and properties 3) and 5) of $\psi$ follows emptiness of the $\Omega$ recession cone. It means that for any nontrivial direction $d \in \mathbb{R}^n$ and any $(\lambda, k) \in \mathbb{R}_{++}^{m+1}$ we have

$$\lim_{t \to \infty} \mathcal{L}(x + td, \lambda, k) = \infty$$

for any $x \in \Omega$. Hence for a given $(\lambda, k) \in \mathbb{R}_{++}^{m+1}$ there exists $\hat{x} \equiv \hat{x}(\lambda, k)$ defined by (100) and $\hat{\lambda} \equiv \hat{\lambda}(\lambda, k)$ defined by (101). Due to 2) of $\psi$ we have $\lambda \in \mathbb{R}_{++}^m \Rightarrow \hat{\lambda} \in \mathbb{R}_{++}^m$, therefore NR method (100)-(101) is well defined.

2) From (100) and (101) follows

$$\nabla_x \mathcal{L}(\hat{x}, \hat{\lambda}, k) = \nabla f(\hat{x}) - \sum_{i=1}^m \lambda_i \psi'(k c_i(\hat{x})) \nabla c_i(\hat{x}) = \nabla_x L(\hat{x}, \hat{\lambda}) = 0,$$

therefore

$$\min_{x \in \mathbb{R}} L(x, \hat{\lambda}) = L(\hat{x}, \hat{\lambda}) = d(\hat{\lambda}).$$

The subdifferential $\partial d(\hat{\lambda})$ contains $-c(\hat{x})$, that is

$$(103) \qquad 0 \in c(\hat{x}) + \partial d(\hat{\lambda}).$$

From (101) follows $\psi'(k c_i(\hat{x})) = \hat{\lambda}_i/\lambda_i, \ i = 1, ..., m.$

Due to 3) of $\psi$ there exists an inverse $\psi'^{-1}$. Using LEID we obtain

(104)     $$c_i(\hat{x}) = k^{-1}\psi'^{-1}(\hat{\lambda}_i/\lambda_i) = k^{-1}\psi^{*'}(\hat{\lambda}_i/\lambda_i)$$

combining (103) and (104) we have

(105)     $$0 \in \partial d(\hat{\lambda}) + k^{-1}\sum_{i=1}^{m}\psi^{*'}\left(\hat{\lambda}_i/\lambda_i\right)e_i.$$

The inclusion (105) is the optimality criteria for $\hat{\lambda}$ to be a solution of problem (102).     $\square$

*Remark* 10. It follows from $1°$ and $2°$ of Assertion 2, that for any $k > 0$ we have $x^* = x(\lambda^*, k)$ and $\lambda^* = \lambda(\lambda^*, k)$, that is $\lambda^* \in \mathbb{R}_+^m$ is a fixed point of the mapping $\lambda \rightarrow \hat{\lambda}(\lambda, k)$.

Along with the class $\Psi$ of transformations $\psi$ we consider a class $\Phi$ of kernels $\varphi = -\psi^*$, with properties induced by properties of $\psi$. We collect them in the following Assertion.

*Assertion* 3. The kernel $\varphi \in \Phi$ are strictly convex on $\mathbb{R}_+$ and possess the following properties on $]0, \infty[$.

  1) $\varphi(s) \geq 0$, $\min_{s \geq 0}\varphi(s) = \varphi(1) = 0$,
  2) $\varphi'(1) = 0$;
  3) $\varphi''(s) > 0$.

Assertion 3 follows from properties 1)-3) of $\psi$ and (11).

The general NR scheme and corresponding methods were introduced in the early 80s (see [46] and references therein). Independently the prox methods with $\varphi$-divergence distance has been studied by M. Teboulle (see [59]). The equivalence of NR and prox methods with $\varphi$- divergence distance was established in [50].

In the following subsection we consider an important particular case of NR - the MBF method.

5.2. **Convergence of the MBF Method and its Dual Equivalent.** For reasons, which will be clear later, we would like to concentrate on the NR method with transformation $\psi(t) = \ln(t + 1)$, which leads to the MBF theory and methods developed in [46] (see also [10], [25], [31], [34], [37], [41], [53] and references therein). The correspondent Lagrangian for the equivalent problem $\mathcal{L} : \mathbb{R}^n \times \mathbb{R}_+^m \times \mathbb{R}_{++} \rightarrow \mathbb{R}$ is defined by formula

$$\mathcal{L}(x, \lambda, k) = f(x) - k^{-1}\sum_{i=1}^{m}\lambda_i\ln(kc_i(x) + 1).$$

For a given $k > 0$ and $\lambda_0 = e = (1, ..., 1) \in \mathbb{R}_{++}^m$ the MBF method generates the following primal-dual sequence $\{x_s, \lambda_s\}_{s=0}^{\infty}$:

$$x_{s+1} : \nabla_x\mathcal{L}(x_{s+1}, \lambda_s, k) =$$

(106)     $$\nabla f(x_{s+1}) - \sum_{i=1}^{m}\lambda_{i,s}(kc_i(x_{s+1}) + 1)^{-1}\nabla c_i(x_{s+1}) = 0$$

(107)     $$\lambda_{s+1} : \lambda_{i,s+1} = \lambda_{i,s}(kc(x_{s+1}) + 1)^{-1}, \ i = 1, ..., m.$$

The Hausdorff distance between two compact sets in $\mathbb{R}_+^m$ will be used later.

Let $X$ and $Y$ be two bounded and closed sets in $\mathbb{R}^n$ and $d(x, y) = \|x - y\|$ is the Euclidean distance between $x \in X, y \in Y$. Then the Hausdorff distance between $X$ and $Y$ is defined as follows

$$d_H(X, Y) := \max\{\max_{x \in X} \min_{y \in Y} d(x, y), \max_{y \in Y} \min_{x \in X} d(x, y)\} =$$

$$\max\{\max_{x \in X} d(x, Y), \max_{y \in Y} d(y, X)\}.$$

For any pair of compact sets $X$ and $Y \subset \mathbb{R}^n$

$$d_H(X, Y) = 0 \Leftrightarrow X = Y.$$

Let $Q \subset \mathbb{R}^m_{++}$ be a compact set, $\hat{Q} = \mathbb{R}^m_{++} \setminus Q$, $S(u, \epsilon) = \{v \in \mathbb{R}^m_+ : \|u - v\| \leq \epsilon\}$ and

$$\partial Q = \{u \in Q | \exists v \in Q : v \in S(u, \epsilon), \exists \hat{v} \in \hat{Q} : \hat{v} \in S(u, \epsilon)\}, \forall \epsilon > 0$$

be the boundary of $Q$.

Let $A \subset B \subset C$ be convex and compact sets in $\mathbb{R}^m_+$. The following inequality follows from the definition of Hausdorff distance.

(108) $$d_H(A, \partial B) < d_H(A, \partial C)$$

Along with the dual sequence $\{\lambda_s\}_{s=0}^{\infty}$ we consider the corresponding convex and bounded level sets $\Lambda_s = \{\lambda \in \mathbb{R}^m_+ : d(\lambda) \geq d(\lambda_s)\}$ and their boundaries $\partial \Lambda_s = \{\lambda \in \Lambda_s : d(\lambda) = d(\lambda_s)\}$.

**Theorem 11.** *Under condition of Theorem 9 for any given $k > 0$ and any $\lambda_0 \in \mathbb{R}^m_{++}$ the MBF method (106)-(107) generates such primal-dual sequence $\{x_s, \lambda_s\}_{s=0}^{\infty}$ that:*

1) $d(\lambda_{s+1}) > d(\lambda_s)$, $s \geq 0$
2) $\lim_{s \to \infty} d(\lambda_s) = d(\lambda^*)$, $\lim_{s \to \infty} f(x_s) = f(x^*)$
3) $\lim_{s \to \infty} d_H(\partial \Lambda_s, L^*) = 0$
4) *there exists a subsequence $\{s_l\}_{l=1}^{\infty}$ such that for $\bar{x}_l = \sum_{s=s_l}^{s_{l+1}} (s_{l+1} - s_l)^{-1} x_s$ we have $\lim_{l \to \infty} \bar{x}_l = \bar{x} \in X^*$, i.e. the primal sequence converges to the primal solution in the ergodic sense.*

*Proof.* 1) It follows from Theorem 9 that method (106)-(107) is well defined and it is equivalent to following proximal point method

$$d(\lambda_{s+1}) - k^{-1} \sum_{i=1}^{m} \lambda_{i,s} \varphi(\lambda_{i,s+1}/\lambda_{i,s}) =$$

(109)

$$\max\{d(u) - k^{-1} \sum_{i=1}^{m} \lambda_{i,s} \varphi(u_i/\lambda_{i,s}) : u \in \mathbb{R}^m_{++}\},$$

where $\varphi = -\psi^* = -\inf_{t > -1}\{st - \ln(t + 1)\} = -\ln s + s - 1$ is the MBF kernel.

The $\varphi$-divergence distance function

$$D(\lambda, u) = \sum_{i=1}^{m} \lambda_i \varphi(u_i/\lambda_i) = \sum_{i=1}^{m} [-\lambda_i \ln u_i/\lambda_i + u_i - \lambda_i],$$

which measures the divergence between two vectors $\lambda$ and $u$ from $\mathbb{R}^m_{++}$ is, in fact, the Kullback-Leibler (KL) distance (see [20], [50], [59]). The MBF kernel $\varphi(s) = -\ln s + s - 1$ is strictly convex on $\mathbb{R}_{++}$ and $\varphi'(1) = 0$, therefore $\min_{s > 0} \varphi(s) = \varphi(1) = 0$, also

a) $D(\lambda, u) > 0, \ \forall \lambda \neq u \in \mathbb{R}_{++}^m$

b) $D(\lambda, u) = 0 \Leftrightarrow \lambda = u$.

From (109) for $u = \lambda_s$ follows

$$(110) \qquad d(\lambda_{s+1}) \geq d(\lambda_s) + k^{-1} \sum_{i=1}^{m} \lambda_{i,s} \varphi(\lambda_{i,s+1}/\lambda_{i,s}).$$

Therefore the sequence $\{d(\lambda_s)\}_{s=0}^{\infty}$ is monotone increasing, unless $\varphi(\lambda_{i,s+1}/\lambda_{i,s}) = 0$ for all $i = 1, ..., m$, but in such case $\lambda_{s+1} = \lambda_s = \lambda^*$. The monotone increasing sequence $\{d(\lambda_s)\}_{s=0}^{\infty}$ is bounded from above by $f(x^*)$, therefore there exists $\lim_{s \to \infty} d(\lambda_s) = \bar{d} \leq f(x^*)$.

2) Our next step is to show that $\bar{d} = f(x^*)$.

From $-c(x_{s+1}) \in \partial d(\lambda_{s+1})$ and concavity of the dual function $d$ follows

$$d(\lambda) - d(\lambda_{s+1}) \leq (-c(x_{s+1}), \lambda - \lambda_{s+1}), \ \forall \lambda \in \mathbb{R}_{++}^m.$$

So for $\lambda = \lambda_s$ we have

$$(111) \qquad d(\lambda_{s+1}) - d(\lambda_s) \geq (c(x_{s+1}), \lambda_s - \lambda_{s+1}).$$

From the update formula (107) follows

$$(112) \qquad (\lambda_{i,s} - \lambda_{i,s+1}) = kc_i(x_{s+1})\lambda_{i,s+1}, \ i = 1, ..., m,$$

therefore from (111) and (112) we have

$$(113) \qquad d(\lambda_{s+1}) - d(\lambda_s) \geq k \sum_{i=1}^{m} c_i^2(x_{s+1})\lambda_{i,s+1}.$$

From Slater condition follows boundedness of $L^*$. Therefore from concavity $d$ follows boundedness of the dual level set

$$\Lambda(\lambda_0) = \{\lambda \in \mathbb{R}_+^m : d(\lambda) \geq d(\lambda_0)\}.$$

It follows from the dual monotonicity (110) that the dual sequence $\{\lambda_s\}_{s=0}^{\infty} \in \Lambda(\lambda_0)$ is bounded.

Therefore there exists $L > 0 : \max_{i,s} \lambda_{i,s} = L$. From (113) follows

$$(114) \qquad d(\lambda_{s+1}) - d(\lambda_s) \geq kL^{-1}(c(x_{s+1}), \lambda_{s+1})^2.$$

By summing up (114) from $s = 1$ to $s = N$ we obtain

$$d(\lambda^*) - d(\lambda_0) \geq d(\lambda_{N+1}) - d(\lambda_0) > kL^{-1} \sum_{s=1}^{N} (\lambda_s, c(x_s))^2,$$

which leads to asymptotic complementarity condition

$$(115) \qquad \lim_{s \to \infty} (\lambda_s, c(x_s)) = 0.$$

On the other hand, from (110) follows

$$(116) \qquad d(\lambda^*) - d(\lambda_0) \geq d(\lambda_N) - d(\lambda_0) \geq k^{-1} \sum_{s=1}^{N} D(\lambda_s, \lambda_{s+1}).$$

Therefore $\lim_{s \to \infty} D(\lambda_s, \lambda_{s+1}) = 0$, which means that divergence (entropy) between two sequential LM vectors asymptotically disappears, that is the dual sequence converges to the fixed point of the map $\lambda \to \hat{\lambda}(\lambda, k)$, which due to Remark 10, is $\lambda^*$.

We need few more steps to prove it. Let us show first that

$$(117) \qquad D(\lambda^*, \lambda_s) > D(\lambda^*, \lambda_{s+1}), \ \forall s \geq 0$$

unless $\lambda_s = \lambda_{s+1} = \lambda^*$.

We assume $x \ln x = 0$ for $x = 0$, then

$$D(\lambda^*, \lambda_s) - D(\lambda^*, \lambda_{s+1}) = \sum_{i=1}^{m} \left( \lambda_i^* \ln \frac{\lambda_{i,s+1}}{\lambda_{i,s}} + \lambda_{i,s} - \lambda_{i,s+1} \right).$$

Invoking the update formula (107) we obtain

$$D(\lambda^*, \lambda_s) - D(\lambda^*, \lambda_{s+1}) = \sum_{i=1}^{m} \lambda_i^* \ln(kc_i(x_{s+1}) + 1)^{-1} + k \sum_{i=1}^{m} \lambda_{i,s+1} c_i(x_{s+1}).$$

Keeping in mind $\ln(1+t)^{-1} = -\ln(1+t) \geq -t$ we have

$$D(\lambda^*, \lambda_s) - D(\lambda^*, \lambda_{s+1}) \geq k \sum_{i=1}^{m} (\lambda_{i,s+1} - \lambda_i^*) c_i(x_{s+1}) =$$

$$(118) \qquad\qquad k(-c(x_{s+1}), \lambda^* - \lambda_{s+1}).$$

From concavity $d$ and $-c(x_{s+1}) \in \partial d(\lambda_{s+1})$ follows

$$(119) \qquad 0 \leq d(\lambda^*) - d(\lambda_{s+1}) \leq (-c(x_{s+1}), \lambda^* - \lambda_{s+1}).$$

Combining (118) and (119) we obtain

$$(120) \qquad D(\lambda^*, \lambda_s) - D(\lambda^*, \lambda_{s+1}) \geq k(d(\lambda^*) - d(\lambda_{s+1})) > 0.$$

Assuming that $d(\lambda^*) - \bar{d} = \rho > 0$ and summing up the last inequality from $s = 0$ to $s = N$ we obtain $D(\lambda^*, \lambda_0) \geq kN\rho$, which is impossible for $N > 0$ large enough.

Therefore $\lim_{s \to \infty} d(\lambda_s) = \bar{d} = d(\lambda^*)$, which together with asymptotic complementarity (115) leads to

$$d(\lambda^*) = \lim_{s \to \infty} d(\lambda_s) = \lim_{s \to \infty} [f(x_s) - (\lambda_s, c(x_s))] =$$

$$(121) \qquad\qquad \lim_{s \to \infty} f(x_s) = f(x^*).$$

3) The dual sequence $\{\lambda_s\}_{s=0}^{\infty}$ is bounded, so it has a converging subsequence $\{\lambda_{s_i}\}_{i=0}^{\infty}$: $\lim_{i \to \infty} \lambda_{s_i} = \bar{\lambda}$. It follows from the dual convergence in value that $\bar{\lambda} = \lambda^* \in L^*$, therefore $\{\lambda \in \mathbb{R}_+^m : d(\lambda) = d(\bar{\lambda})\} = L^*$.

From (110) follows $L^* \subset ... \subset \Lambda_{s+1} \subset \Lambda_s \subset ... \subset \Lambda_0$, therefore from (108) we obtain a monotone decreasing sequence $\{d_H(\partial \Lambda_s, L^*)\}_{s=0}^{\infty}$, which has a limit, that is

$$\lim_{s \to \infty} d_H(\partial \Lambda_s, L^*) = \rho \geq 0,$$

but $\rho > 0$ is impossible due to the continuity of the dual function and the convergence of the dual sequence in value.

4) Let us consider the indices subset $I_+ = \{i : \bar{\lambda}_i > 0\}$, then from (115) we have $\lim_{s \to \infty} c_i(x_s) = c_i(\bar{x}) = 0, \quad i \in I_+$. Now we consider the indices subset $I_0 = \{i : \bar{\lambda}_i = 0\}$.

There exists a subsequence $\{\lambda_{s_l}\}_{l=1}^{\infty}$ that $\lambda_{i,s_{l+1}} \leq 0.5\lambda_{i,s_l}, i \in I_0$.

Using again the update formula (107) we obtain

$$\lambda_{s_{l+1}} \prod_{s=s_l}^{s_{l+1}} (kc_i(x_s) + 1) = \lambda_{i,s_l} \geq 2\lambda_{s_{l+1}}, \ i \in I_0.$$

Invoking the arithmetic-geometric means inequality we have

$$\frac{1}{s_{l+1} - s_l} \sum_{s=s_l}^{s_{l+1}} (kc_i(x_s) + 1) \geq \left( \prod_{s=s_l+1}^{s_{l+1}} (kc_i(x_s) + 1) \right)^{1/(s_{l+1} - s_l)} \geq 2^{(1/s_{l+1} - s_l)} > 1.$$

Therefore

$$\frac{k}{(s_{l+1} - s_l)} \sum_{s=s_l}^{s_{l+1}} c_i(x_s) > 0 \ i \in I_0.$$

From concavity $c_i$ we obtain

$$(122) \quad c_i(\bar{x}_{l+1}) = c_i \left( \sum_{s=s_l+1}^{s_{l+1}} \frac{1}{s_{l+1} - s_l} x_s \right) \geq \frac{1}{s_{l+1} - s_l} \sum_{s=s_l+1}^{s_{l+1}} c_i(x_s) > 0, \ i \in I_0.$$

On the other hand, from convexity of $f$ we have

$$(123) \quad f(\bar{x}_{l+1}) \leq \frac{1}{s_{l+1} - s_l} \sum_{s=s_l+1}^{s_{l+1}} f(x_s).$$

Without loosing generality we can assume that $\lim_{l\to\infty} \bar{x}_l = \bar{x} \in \Omega$. It follows from (121) that

$$f(\bar{x}) = \lim_{l\to\infty} f(\bar{x}_l) \leq \lim_{s\to\infty} f(x_s) = \lim_{s\to\infty} d(\lambda_s) = d(\lambda^*) = f(x^*).$$

Thus $f(\bar{x}) = f(x^*) = d(\lambda^*) = d(\bar{\lambda})$ and $\bar{x} = x^*$, $\bar{\lambda} = \lambda^*$. The proof of Theorem 11 is completed. $\qquad\square$

We conclude the section with few remarks.

*Remark* 12. Each $\psi \in \Psi$ leads to a particular NR method for solving (51) as well as to an interior prox method for solving the dual problem (53). In this regard NR approach is source of methods for solving (53), which arises in a number of application such as non-negative least square, statistical learning theory, image space reconstruction, maximum likelihood estimation in emission tomography (see [17], [20], [62] and references therein).

*Remark* 13. The MBF method leads to the multiplicative method (107) for the dual problem. If the dual function $d$ has a gradient, then $\nabla d(\lambda_{s+1}) = -c(x_{s+1})$. Formulas (107) can be rewritten as follows

$$(124) \quad \lambda_{i,s+1} - \lambda_{i,s} = k\lambda_{i,s+1}[\nabla d(\lambda_{s+1})], \ i = 1, ..., m,$$

which is, in fact, implicit Euler method for the following system of ordinary differential equations

$$(125) \quad \frac{d\lambda}{dt} = k\lambda\nabla d(\lambda), \ \lambda(0) = \lambda_0.$$

Therefore the dual MBF method (124) is called (see (1.7) in [20]) implicit multiplicative algorithm.

The explicit multiplicative algorithm (see (1.8) in [20]) is given by the following formula

$$\lambda_{i,s+1} = \lambda_{i,s}(1 - k[\nabla d(\lambda_s)]_i)^{-1}, \ i = 1, ..., m. \tag{126}$$

It has been used by Eggermond [20] for solving non-negative least square, by Daube-Witherspoon and Muchlehner [17] for image space reconstruction (ISRA) and by Shepp and Vardi in their EM method for finding maximum likelihood estimation in emission tomography [62].

*Remark* 14. Under the standard second order sufficient optimality condition there exists $k_0 > 0$ that for $k \geq k_0$ the MBF method (106)-(107) converges with linear rate

$$\|x_{s+1} - x^*\| \leq \frac{c}{k}\|\lambda_s - \lambda^*\|; \ \|\lambda_{s+1} - \lambda^*\| \leq \frac{c}{k}\|\lambda_s - \lambda^*\|$$

and $c > 0$ is independent on $k \geq k_0$. By increasing $k$ from step to step one obtains superlinear convergence rate (see [46]).

## 6. Lagrangian Transformation and Interior ellipsoid methods

The Lagrangian transformation (LT) scheme employs a class $\psi \in \Psi$ of smooth strictly concave, monotone increasing functions to transform terms of the Classical Lagrangian associated with constraints. The transformation is scaled by a positive scaling parameter.

Finding a primal minimizer of the transformed Lagrangian following by the Lagrange multipliers update leads to a new class of multipliers methods.

The LT methods are equivalent to proximal point methods with Bregman or Bregman type distance function for the dual problem. The kernel of the correspondent distance is $\varphi = -\psi^*$.

Each dual prox, in turn, is equivalent to an interior ellipsoid methods. In case of the MBF transformation $\psi(t) = \ln(t+1)$ the dual prox is based on Bregman distance $B(u, v) = \sum_{i=1}^{m}(-\ln(u_i/v_i) + u_i/v_i - 1)$ with MBF kernel $\varphi = -\psi^* = -\ln s + s - 1$, which is SC function. Therefore the interior ellipsoids are Dikin's ellipsoids (see [18], [39], [42], [43], [49] ).

Application of LT with MBF transformation for LP leads to Dikin's affine scaling type method for the dual LP.

6.1. **Lagrangian Transformation.** We consider a class $\Psi$ of twice continuous differentiable functions $\psi : R \to R$ with the following properties

1) $\psi(0) = 0$
2) a) $\psi'(t) > 0$,    b) $\psi'(0) = 1$,    $\psi'(t) \leq at^{-1}, a > 0, t > 0$
3) $-m_0^{-1} \leq \psi''(t) < 0$,    $\forall t \in ]-\infty, \infty[$
4) $\psi''(t) \leq -M^{-1}$,    $\forall t \in ]-\infty, 0[$,    $0 < m_0 < M < \infty$.

For a given $\psi \in \Psi$ and $k > 0$, the LT $\mathcal{L} : \mathbb{R}^n \times \mathbb{R}_+^m \times \mathbb{R}_{++} \to \mathbb{R}$ is defined by the following formula

$$\mathcal{L}(x, \lambda, k) := f(x) - k^{-1}\sum_{i=1}^{m}\psi(k\lambda_i c_i(x)). \tag{127}$$

It follows from 2a) and 3), convexity $f$, concavity $c_i, i = 1, ..., m$ that for any given $\lambda \in \mathbb{R}_{++}^m$ and any $k > 0$ the LT is convex in $x$.

6.2. **Primal Transformations and Dual Kernels.** The well known transformations

- exponential [7], [40], [61] $\hat{\psi}_1(t) = 1 - e^{-t}$;
- logarithmic MBF [46] $\hat{\psi}_2(t) = \ln(t+1)$;
- hyperbolic MBF [46] $\hat{\psi}_3(t) = t/(t+1)$;
- log-sigmoid [48] $\hat{\psi}_4(t) = 2(\ln 2 + t - \ln(1 + e^t))$;
- Chen-Harker-Kanzow-Smale [48] (CHKS) $\hat{\psi}_5(t) = t - \sqrt{t^2 + 4\eta} + 2\sqrt{\eta}$, $\eta > 0$, unfortunately, do not belong to $\Psi$.

The transformations $\hat{\psi}_1, \hat{\psi}_2, \hat{\psi}_3$ do not satisfy 3) ($m_0 = 0$), while transformations $\hat{\psi}_4$ and $\hat{\psi}_5$ do not satisfy 4) ($M = \infty$). A slight modification of $\hat{\psi}_i, i = 1, \ldots, 5$, however, leads to $\psi_i \in \Psi$ (see [6]).

Let $-1 < \tau < 0$, we will use later the following truncated transformations $\psi_i : \mathbb{R} \to \mathbb{R}$ are defined as follows

$$(128) \qquad \psi_i(t) := \begin{cases} \hat{\psi}_i(t), \infty > t \geq \tau \\ q_i(t), -\infty < t \leq \tau, \end{cases}$$

where $q_i(t) = a_i t^2 + b_i t + c_i$ and $a_i = 0.5\hat{\psi}_i''(\tau)$, $b_i = \hat{\psi}_i'(\tau) - \tau\hat{\psi}''(\tau)$, $c_i = \hat{\psi}_i'(\tau) - \tau\hat{\psi}_i'(\tau) + 0.5\tau^2\hat{\psi}_i''(\tau)$.

It is easy to check that for truncated transformations $\psi_i$, $i = 1, \ldots, 5$ the properties 1)-4) hold, that is $\psi_i \in \Psi$.

In the future along with transformations $\psi_i \in \Psi$ their conjugate

$$(129) \qquad \psi_i^*(s) := \begin{cases} \hat{\psi}_i^*(s), & s \leq \hat{\psi}_i'(\tau) \\ q_i^*(s) = (4a_i)^{-1}(s - b_i)^2 - c_i, & s \geq \hat{\psi}_i'(\tau), i = 1, \ldots, 5, \end{cases}$$

will play an important role, where $\hat{\psi}^*{}_i(s) = \inf_t \{st - \hat{\psi}_i(t)\}$ is the LET of $\hat{\psi}_i$.

With the class of primal transformations $\Psi$ we associate the class of dual kernels

$$\varphi \in \Phi = \{\varphi = -\psi^* : \psi \in \Psi\}.$$

Using properties 2) and 4) one can find $0 < \theta_i < 1$ that

$$\hat{\psi}_i'(\tau) - \hat{\psi}_i'(0) = -\hat{\psi}_i''(\tau\theta_i)(-\tau) \geq -\tau M^{-1}, i = 1, \ldots, 5$$

or

$$\hat{\psi}_i'(\tau) \geq 1 - \tau M^{-1} = 1 + |\tau| M^{-1}.$$

Therefore from (129) for any $0 < s \leq 1 + |\tau| M^{-1}$ we have

$$(130) \qquad \varphi_i(s) = \hat{\varphi}_i(s) = -\hat{\psi}_i^*(s) = \inf_t \{st - \hat{\psi}_i(t)\},$$

where kernels

- exponential $\hat{\varphi}_1(s) = s \ln s - s + 1, \hat{\varphi}_1(0) = 1$;
- logarithmic MBF $\hat{\varphi}_2(s) = -\ln s + s - 1$;
- hyperbolic MBF $\hat{\varphi}_3(s) = -2\sqrt{s} + s + 1, \hat{\varphi}_3(0) = 1$;
- Fermi-Dirac $\hat{\varphi}_4(s) = (2 - s)\ln(2 - s) + s \ln s, \hat{\varphi}_4(0) = 2\ln 2$;
- CMKS $\hat{\varphi}_5(s) = -2\sqrt{\eta}(\sqrt{(2 - s)s} - 1), \hat{\varphi}_5(0) = 2\sqrt{\eta}$

  are infinitely differentiable on $]0, 1 + |\tau| M^{-1}[$.

To simplify the notations we omit indices of $\psi$ and $\varphi$.

The properties of kernels $\varphi \in \Phi$ induced by 1)–4) can be established by using (11).

We collect them in the following Assertion

*Assertion* 4. The kernels $\varphi \in \Phi$ are strictly convex on $\mathbb{R}_+^m$, twice continuously differentiable and possess the following properties

   1) $\varphi(s) \geq 0, \ \forall s \in ]0, \infty[$ and $\min\limits_{s \geq 0} \varphi(s) = \varphi(1) = 0$;

   2) a) $\lim\limits_{s \to 0^+} \varphi'(s) = -\infty,$    b) $\varphi'(s)$ is monotone increasing and
         c) $\varphi'(1) = 0$;

   3) a) $\varphi''(s) \geq m_0 > 0, \quad \forall s \in ]0, \infty[,$    b) $\varphi''(s) \leq M < \infty, \quad \forall s \in [1, \infty[$.

Let $Q \subset \mathbb{R}^m$ be an open convex set, $\hat{Q}$ is the closure of $Q$ and $\varphi : \hat{Q} \to \mathbb{R}$ be a strictly convex closed function on $\hat{Q}$ and continuously differentiable on $Q$, then the Bregman distance $\mathbb{B}_\varphi : \hat{Q} \times Q \to R_+$ induced by $\varphi$ is defined as follows(see [8]),

$$(131) \qquad \mathbb{B}_\varphi(x, y) = \varphi(x) - \varphi(y) - (\nabla\varphi(y), x - y).$$

Let $\varphi \in \Phi$, then $B_\varphi : \mathbb{R}_+^m \times \mathbb{R}_{++}^m \to \mathbb{R}_+$, defined by

$$B_\varphi(u, v) := \sum_{i=1}^{m} \varphi(u_i/v_i),$$

we call Bregman type distance induced by kernel $\varphi$. Due to $\varphi(1) = \varphi'(1) = 0$ for any $\varphi \in \Phi$, we have

$$(132) \qquad \varphi(t) = \varphi(t) - \varphi(1) - \varphi'(1)(t - 1),$$

which means that $\varphi(t) : \mathbb{R}_{++} \to \mathbb{R}_{++}$ is Bregman distance between $t > 0$ and 1.

By taking $t_i = \frac{u_i}{v_i}$ from (132) we obtain

$$(133) \qquad B_\varphi(u, v) = B_\varphi(u, v) - B_\varphi(v, v) - (\nabla_u B_\varphi(v, v), u - v),$$

which justifies the definition of the Bregman type distance.

For the MBF kernel $\varphi_2(s) = -\ln s + s - 1$ we obtain the Bregman distance,

$$(134) \qquad \begin{aligned} \mathbb{B}_2(u, v) &= \sum_{i=1}^{m} \varphi_2(u_i/v_i) = \sum_{i=1}^{m}(-\ln u_i/v_i + u_i/v_i - 1) = \\ &\sum_{i=1}^{m}[-\ln u_i + \ln v_i + (u_i - v_i)/v_i], \end{aligned}$$

which is induced by the standard log-barrier function $F(t) = -\sum_{i=1}^{m} \ln t_i$.

After Bregman's introduction his function in the 60s (see [8]) the prox method with Bregman distance has been widely studied (see [9], [11], [13], [15], [19], [39], [48]-[50] and reference therein).

From the definition of $\mathbb{B}_2(u, v)$ follows

$$\nabla_u \mathbb{B}_2(u, v) = \nabla F(u) - \nabla F(v).$$

For $u \in \hat{Q}$, $v \in Q$ and $w \in Q$ the following three point identity established by Chen and Teboulle in [15] is an important element in the analysis of prox methods with Bregman distance

$$(135) \qquad \mathbb{B}_2(u, v) - \mathbb{B}_2(u, w) - \mathbb{B}_2(w, v) = (\nabla F(v) - \nabla F(w), w - u).$$

The properties of Bregman's type distance functions we collect in the following Assertion.

*Assertion* 5. The Bregman type distance satisfies the following properties:

1) $B_\varphi(u, v) \geq 0, \forall u \in \mathbb{R}_+^m, v \in \mathbb{R}_{++}^m, \ B_\varphi(u, v) = 0 \Leftrightarrow u = v, \ \forall v, u \in \mathbb{R}_{++}^m;$ $B_\varphi(u, v) > 0$ for any $u \neq v$

2) $B_\varphi(u, v) \geq \frac{1}{2} m_0 \sum_{i=1}^m (\frac{u_i}{v_i} - 1)^2, \forall u_i \in \mathbb{R}_+^m, v_i \in \mathbb{R}_{++}^m;$

3) $B_\varphi(u, v) \leq \frac{1}{2} M \sum_{i=1}^m (\frac{u_i}{v_i} - 1)^2, \forall u \in \mathbb{R}_+^m, u \geq v > 0;$

4) for any fixed $v \in \mathbb{R}_{++}^m$ the gradient $\nabla_u B_\varphi(u, v)$ is a barrier function of $u \in \mathbb{R}_{++}^m, \ i.e.$

$$\lim_{u_i \to 0_+} \frac{\partial}{\partial u_i} B_\varphi(u, v) = -\infty, i = 1, \dots, m.$$

The properties 1)–4) directly following from the properties of kernels $\varphi \in \Phi$ given in Assertion 4.

### 6.3. Primal LT and Dual Prox Methods.
Let $\psi \in \Psi$, $\lambda_0 \in \mathbb{R}_{++}^m$ and $k > 0$ are given. The LT method generates a primal–dual sequence $\{x_s, \lambda_s\}_{s=1}^\infty$ by formulas

$$(136) \qquad x_{s+1} : \nabla_x \mathcal{L}(x_{s+1}, \lambda_s, k) = 0$$

$$(137) \qquad \lambda_{i,s+1} = \lambda_{i,s} \psi'(k\lambda_{i,s} c_i(x_{s+1})), i = 1, \dots, m.$$

**Theorem 15.** *If conditions A and B hold and $f$, $c_i$, $i = 1, ..., m$ continuously differentiable then:*

1) *the LT method (136)-(137) is well defined and it is equivalent to the following interior proximal point method*

$$(138) \qquad \lambda_{s+1} = \mathrm{argmax}\{d(\lambda) - k^{-1} B_\varphi(\lambda, \lambda_s) | \lambda \in \mathbb{R}_{++}^m\},$$

*where*

$$B_\varphi(u, v) := \sum_{i=1}^m \varphi(u_i/v_i)$$

*and $\varphi = -\psi^*.$*

2) *for all $i = 1, ..., m$ we have*

$$(139) \qquad \lim_{s \to \infty} (\lambda_{i,s+1}/\lambda_{i,s}) = 1.$$

*Proof.*     1) From assumptions A, convexity of $f$, concavity of $c_i, i = 1, \dots, m$ and property 4) of $\psi \in \Psi$ for any $\lambda_s \in \mathbb{R}_{++}^m$ and $k > 0$ follows boundedness of the level set $\{x : \mathcal{L}(x, \lambda_s, k) \leq \mathcal{L}(x_s, \lambda_s, k)\}$. Therefore, the minimizer $x_s$ exists for any $s \geq 1$. It follows from 2 a) of $\psi \in \Psi$ and (137) that $\lambda_s \in \mathbb{R}_{++}^m \Rightarrow \lambda_{s+1} \in \mathbb{R}_{++}^m$. Therefore the LT method (136)– (137) is well defined. From (136) follows

$$\nabla_x \mathcal{L}(x_{s+1}, \lambda_s, k) =$$

$$(140) \qquad \nabla f(x_{s+1}) - \sum_{i=1}^m \lambda_{i,s} \psi'(k\lambda_{i,s} c_i(x_{s+1})) \nabla c_i(x_{s+1}) = 0.$$

From (136) and (137) we obtain

$$\nabla_x \mathcal{L}(x_{s+1}, \lambda_s, k) = \nabla f(x_{s+1}) - \sum_{i=1}^m \lambda_{i,s+1} \nabla c_i(x_{s+1}) = \nabla_x L(x_{s+1}, \lambda_{s+1}) = 0,$$

therefore

$$d(\lambda_{s+1}) = L(x_{s+1}, \lambda_{s+1}) = \min\{L(x, \lambda_{s+1}) | x \in \mathbb{R}^n\}.$$

From (137) we get

$$\psi'(k\lambda_{i,s} c_i(x_{s+1})) = \lambda_{i,s+1}/\lambda_{i,s}, i = 1, \dots, m.$$

In view of property 3) for any $\psi \in \Psi$ there exists an inverse $\psi'^{-1}$, therefore

(141) $$c_i(x_{s+1}) = k^{-1}(\lambda_{i,s})^{-1} \psi'^{-1}(\lambda_{i,s+1}/\lambda_{i,s}), i = 1, \dots, m.$$

Using LEID $\psi'^{-1} = \psi^{*\prime}$ we obtain

(142) $$c_i(x_{s+1}) = k^{-1}(\lambda_{i,s})^{-1} \psi^{*\prime}(\lambda_{i,s+1}/\lambda_{i,s}), \ \ i = 1, \dots, m.$$

Keeping in mind

$$-c(\lambda_{s+1}) \in \partial d(\lambda_{s+1})$$

and $\varphi = -\psi^*$ we have

$$0 \in \partial d(\lambda_{s+1}) - k^{-1} \sum_{i=1}^{m} (\lambda_{i,s})^{-1} \varphi'(\lambda_{i,s+1}/\lambda_{i,s}) e_i.$$

The last inclusion is the optimality criteria for $\lambda_{s+1} \in \mathbb{R}_{++}^m$ to be the solution of the problem (138). Thus, the LT method (136)-(137) is equivalent to the interior proximal point method (138).

2) From 1) of Assertion 5 and (138) follows

(143) $$d(\lambda_{s+1}) \geq k^{-1} B_\varphi(\lambda_{s+1}, \lambda_s) + d(\lambda_s) > d(\lambda_s), \ \forall s > 0.$$

Summing up last inequality from $s = 0$ to $s = N$, we obtain

$$d(\lambda^*) - d(\lambda_0) \geq d(\lambda_{N+1}) - d(\lambda_0) > k^{-1} \sum_{s=0}^{N} B_\varphi(\lambda_{s+1}, \lambda_s),$$

therefore

(144) $$\lim_{s \to \infty} B(\lambda_{s+1}, \lambda_s) = \lim_{s \to \infty} \sum_{i=1}^{m} \varphi(\lambda_{i,s+1}/\lambda_{i,s}) = 0.$$

From (144) and 2) of Assertion 5 follows

(145) $$\lim_{s \to \infty} \lambda_{i,s+1}/\lambda_{i,s} = 1, \ i = 1, \dots, m.$$

$\square$

*Remark* 16. From (130) and (145) follows that for $s \geq s_0 > 0$ the Bregman type distance functions $B_{\varphi_i}$ used in (138) are based on kernels $\varphi_i$, which correspond to the original transformations $\hat{\psi}_i$.

The following Theorem establishes the equivalence of LT multipliers method and interior ellipsoid methods (IEMs) for the dual problem.

**Theorem 17.** *It conditions of Theorem 15 are satisfied then:*

1) *for a given $\varphi \in \Phi$ there exists a diagonal matrix $H_\varphi = diag(h_\varphi^i)_{i=1}^m$ with $h_\varphi^i > 0, i = 1, \dots, m$ that $B_\varphi(u, v) = \frac{1}{2} \|u - v\|_{H_\varphi}^2$, where $\|w\|_{H_\varphi}^2 = w^T H_\varphi w$;*

2) *The Interior Prox method (138) is equivalent to an interior quadratic prox (IQP) in the rescaled from step to step dual space, i.e.*

(146)
$$\lambda_{s+1} = \arg\max\{d(\lambda) - \frac{1}{2k}\|\lambda - \lambda_s\|_{H^s_\varphi}^2 | \lambda \in \mathbb{R}^m_+\},$$

*where $H^s_\varphi = diag(h^{i,s}_\varphi) = diag(2\varphi''(1 + \theta^s_i(\lambda_{i,s+1}/\lambda_{i,s} - 1))(\lambda_{i,s})^{-2})$ and $0 < \theta^s_i < 1$;*

3) *The IQP is equivalent to an interior ellipsoid method (IEM) for the dual problem;*

4) *There exists a converging to zero sequence $\{r_s > 0\}_{s=0}^{\infty}$ and step $s_0 > 0$ such that, for $\forall s \geq s_0$, the LT method (136)– (137) with truncated MBF transformation $\psi_2(t)$ is equivalent to the following IEM for the dual problem*

(147)
$$\lambda_{s+1} = \arg\max\{d(\lambda)|\lambda \in E(\lambda_s, r_s)\},$$

*where $H_s = diag(\lambda_{i,s})_{i=1}^{m}$ and*

$$E(\lambda_s, r_s) = \{\lambda : (\lambda - \lambda_s)^T H_s^{-2}(\lambda - \lambda_s) \leq r_s^2\}$$

*is Dikin's ellipsoid associated with the standard log–barrier function $F(\lambda) = -\sum_{i=1}^{m} \ln \lambda_i$ for the dual feasible set $\mathbb{R}^m_+$.*

*Proof.*     1) It follows from $\varphi(1) = \varphi'(1) = 0$ that

(148)
$$B_\varphi(u, v) = \frac{1}{2} \sum_{i=1}^{m} \varphi''(1 + \theta_i(\frac{u_i}{v_i} - 1))(\frac{u_i}{v_i} - 1)^2,$$

where $0 < \theta_i < 1$, $i = 1, \ldots, m$.

Due to 3a) from Assertion 4, we have $\varphi''(1 + \theta_i(\frac{u_i}{v_i} - 1)) \geq m_0 > 0$, and due to property 2a) of $\psi \in \Psi$, we have $v \in \mathbb{R}^m_{++}$, therefore

$$h^i_\varphi = 2\varphi''(1 + \theta_i(\frac{u_i}{v_i} - 1))v_i^{-2} > 0, i = 1, \ldots, m.$$

We consider the diagonal matrix $H_\varphi = diag(h^i_\varphi)_{i=1}^{m}$, then from (148) we have

(149)
$$B_\varphi(u, v) = \frac{1}{2}\|u - v\|_{H_\varphi}^2.$$

2) By taking $u = \lambda$, $v = \lambda_s$ and $H_\varphi = H^s_\varphi$ from (138) and (149) we obtain (146)

3) Let's consider the optimality criteria for the problem (146). Keeping in mind $\lambda_{s+1} \in \mathbb{R}^m_{++}$ we conclude that $\lambda_{s+1}$ is an unconstrained maximizer in (146). Therefore one can find $g_{s+1} \in \partial d(\lambda_{s+1})$ that

(150)
$$g_{s+1} - k^{-1}H^s_\varphi(\lambda_{s+1} - \lambda_s) = 0.$$

Let $r_s = \|\lambda_{s+1} - \lambda_s\|_{H^s_\varphi}$, we consider an ellipsoid

$$E_\varphi(\lambda_s, r_s) = \{\lambda : (\lambda - \lambda_s)^T H^s_\varphi(\lambda - \lambda_s) \leq r_s^2\}$$

with center $\lambda_s \in \mathbb{R}^m_{++}$ and radius $r_s$. It follows from 4) of Assertion 5 that $E(\lambda_s, r_s)$ is an interior ellipsoid in $\mathbb{R}^m_{++}$, i.e. $E_\varphi(\lambda_s, r_s) \subset \mathbb{R}^m_{++}$.

Moreover $\lambda_{s+1} \in \partial E_\varphi(\lambda_s, r_s) = \{\lambda : (\lambda - \lambda_s)^T H^s_\varphi(\lambda - \lambda_s) = r_s^2\}$, therefore (150) is the optimality condition for the following optimization problem

(151)
$$d(\lambda_{s+1}) = \max\{d(\lambda)|\lambda \in E_\varphi(\lambda_s, r_s)\}$$

and $(2k)^{-1}$ is the optimal Lagrange multiplier for the only constraint in (151).

Thus, the Interior Prox method (138) is equivalent to the IEM (151).

4) Let us consider the LT method (136)-(137) with truncated MBF transformation.

From (139) follows that for $s \geq s_0$ only Bregman distance

$$\mathbb{B}_2(\lambda, \lambda_s) = \sum_{i=1}^{m} (-ln\frac{\lambda_i}{\lambda_i^s} + \frac{\lambda_i}{\lambda_i^s} - 1)$$

is used in the LT method (136)-(137). Then

$$\nabla_{\lambda\lambda}^2 \mathbb{B}_2(\lambda, \lambda_s)|_{\lambda=\lambda_s} = H_s^{-2} = (I \cdot \lambda_s)^{-2}.$$

In view of $\mathbb{B}_2(\lambda_s, \lambda_s) = 0$ and $\nabla_\lambda \mathbb{B}_2(\lambda_s, \lambda_s) = 0^m$, we obtain

$$\mathbb{B}_2(\lambda, \lambda_s) = \frac{1}{2}(\lambda - \lambda_s)^T H_s^{-2}(\lambda - \lambda_s) + o(\|\lambda - \lambda_s\|^2) =$$

$$= Q(\lambda, \lambda_s) + o(\|\lambda - \lambda_s\|^2).$$

It follows from (139) that for a any $s \geq s_0$ the term $o(\|\lambda_{s+1} - \lambda_s\|^2)$ can be ignored. Then the optimality criteria (150) can be rewritten as follows

$$g_{s+1} - k^{-1} H_s^{-2}(\lambda_{s+1} - \lambda_s) = 0.$$

Therefore

$$d(\lambda_{s+1}) = \max\{d(\lambda) | \lambda \in E(\lambda_s, r_s)\},$$

where $r_s^2 = Q(\lambda_{s+1}, \lambda_s)$ and

$$E(\lambda_s, r_s) = \{\lambda : (\lambda - \lambda_s)H_s^{-2}(\lambda - \lambda_s) = r_s^2\}$$

is Dikin's ellipsoid. The proof is completed                          □

□

The results of Theorem 17 were used in [49] for proving convergence LT method (136)-(137) and its dual equivalent (138) for Bregman type distance function.

Now we consider the LT method with truncated MBF transformation $\psi_2$.

It follows from (130) and (139) that for $s \geq s_0$ only original transformation $\psi_2(t) = \ln(t+1)$ is used in LT method (136)-(137), therefore only Bregman distance $\mathbb{B}_2(u, v) = \sum_{i=1}^{m}(-\ln(u_i/v_i) + u_i/v_i - 1)$ is used in the prox method (138).

In other words, for a given $k > 0$ the primal-dual sequence $\{x_s, \lambda_s\}_{s=s_0}^{\infty}$ is generated by the following formulas

$$x_{s+1} : \nabla_k \mathcal{L}(x_{s+1}, \lambda_s, k) =$$

(152)

$$\nabla f(x_{s+1}) - \sum_{i=1}^{m} \lambda_{i,s}(1 + k\lambda_{i,s}c_i(x_{s+1}))^{-1}\nabla c_i(x_{s+1}) = 0$$

(153)          $$\lambda_{s+1} : \lambda_{i,s+1} = \lambda_{i,s}(1 + k\lambda_{i,s}c_i(x_{s+1}))^{-1}, \ i = 1, ..., m.$$

The method (152)-(153) Matioti and Gonzaga called $M^2BF$ (see [39]).

**Theorem 18.** *Under condition of Theorem 15 the $M^2BF$ method generates such primal-dual sequence that:*

1) $d(\lambda_{s+1}) > d(\lambda_s), \ s \geq s_0$

2) *a)* $\lim_{s \to \infty} d(\lambda_s) = d(\lambda^*);$ *b)* $\lim_{s \to \infty} f(x_s) = f(x^*)$ *and*

$$c) \quad \lim_{s \to \infty} d_H(\partial\Lambda_s, L^*) = 0$$

3) *there is a subsequence* $\{s_l\}_{l=1}^\infty$ *that for* $\bar{\lambda}_{i,s} = \lambda_{i,s} \left(\sum_{s=s_l}^{s_{l+1}} \lambda_{i,s}\right)^{-1}$ *the sequence* $\{\bar{x}_{l+1} = \sum_{s=s_l}^{s_{l+1}} \bar{\lambda}_{i,s} x_s\}_{l=0}^\infty$ *converges and* $\lim_{l\to\infty} \bar{x}_l = \bar{x} \in X^*$.

*Proof.* 1) From Theorem 17 follows that LT (152)-(153) is equivalent to the prox method (138) with Bregman distance. From (138) with $\lambda = \lambda_s$ we obtain

(154) $$d(\lambda_{s+1}) \geq d(\lambda_s) + k^{-1} \sum_{i=1}^m (-\ln(\lambda_{i,s+1}/\lambda_{i,s}) + \lambda_{i,s+1}/\lambda_{i,s} - 1).$$

The Bregman distance is strictly convex in $u$, therefore from (154) follows $d(\lambda_{s+1}) > d(\lambda_s)$ unless $\lambda_{s+1} = \lambda_s \in \mathbb{R}_{++}^m$, then $c_i(x_{s+1}) = 0$, $i = 1, .., m$ and $(x_{s+1}, \lambda_{s+1}) = (x^*, \lambda^*)$ is a KKT pair.

2) The monotone increasing sequence $\{d(\lambda_s)\}_{s=s_0}^\infty$ is bounded from above by $f(x^*)$, therefore there exists $\bar{d} = \lim_{s\to\infty} d(\lambda_s) \leq d(\lambda^*) = f(x^*)$.

The first step is to show that $\bar{d} = d(\lambda^*)$.

Using $\nabla_u \mathbb{B}_2(v, w) = \nabla F(v) - \nabla F(w)$ for $v = \lambda_s$ and $w = \lambda_{s+1}$ we obtain

$$\nabla_\lambda \mathbb{B}_2(\lambda, \lambda_{s+1})_{/\lambda=\lambda_s} = \nabla\varphi(\lambda_s) - \nabla\varphi(\lambda_{s+1}) = \left(-\sum_{i=1}^m \lambda_{i,s}^{-1} e_i + \sum_{i=1}^m \lambda_{i,s+1}^{-1} e_i\right).$$

From the three point identity (135) with $u = \lambda^*$, $v = \lambda_s$, $w = \lambda_{s+1}$ follows

(155)
$$\mathbb{B}_2(\lambda^*, \lambda_s) - \mathbb{B}_2(\lambda^*, \lambda_{s+1}) - \mathbb{B}_2(\lambda_{s+1}, \lambda_s) =$$
$$(\nabla\varphi(\lambda_s) - \nabla\varphi(\lambda_{s+1}), \lambda_{s+1} - \lambda^*) =$$
$$\sum_{i=1}^m (-\lambda_{i,s}^{-1} + \lambda_{i,s+1}^{-1})(\lambda_{i,s+1} - \lambda_i^*).$$

From the update formula (153) follows

$$kc_i(x_{s+1}) = -\lambda_{i,s}^{-1} + \lambda_{i,s+1}^{-1}, \ i = 1, ..., m.$$

Therefore, keeping in mind, $\mathbb{B}_2(\lambda_s, \lambda_{s+1}) \geq 0$ we can rewrite (155) as follows

$$\mathbb{B}_2(\lambda^*, \lambda_s) - \mathbb{B}_2(\lambda^*, \lambda_{s+1}) \geq k(c(x_{s+1}), \lambda_{s+1} - \lambda^*).$$

From $-c(x_{s+1}) \in \partial d(\lambda_{s+1})$ we obtain

(156) $$d(\lambda) - d(\lambda_{s+1}) \leq (-c(x_{s+1}), \lambda - \lambda_{s+1}), \forall \lambda \in \mathbb{R}_+^m.$$

For $\lambda = \lambda^*$ from (156) we get

$$(c(x_{s+1}), \lambda_{s+1} - \lambda^*) \geq d(\lambda^*) - d(\lambda_{s+1}).$$

Hence,

(157) $$\mathbb{B}_2(\lambda^*, \lambda_s) - \mathbb{B}_2(\lambda^*, \lambda_{s+1}) \geq k(d(\lambda^*) - d(\lambda_{s+1})).$$

Assuming $\lim_{s\to\infty} d(\lambda_s) = \bar{d} < d(\lambda^*)$ we have $d(\lambda^*) - d(\lambda_s) \geq \rho > 0, \forall s \geq s_0$. Summing up (157) from $s = s_0$ to $s = N$ we obtain

$$\mathbb{B}_2(\lambda^*, \lambda_{s_0}) - k(N - s_0)\rho \geq \mathbb{B}_2(\lambda^*, \lambda_{N+1}),$$

which is impossible for large $N$. Therefore

(158) $$\lim_{s\to\infty} d(\lambda_s) = d(\lambda^*).$$

From (156) with $\lambda = \lambda_s$ we obtain

$$d(\lambda_s) - d(\lambda_{s+1}) \leq (-c(x_{s+1}), \lambda_s - \lambda_{s+1}).$$

Using the update formula (153) from last inequality we obtain

$$d(\lambda_{s+1}) - d(\lambda_s)) \geq (c(x_{s+1}), \lambda_s - \lambda_{s+1}) =$$

$$(159) \qquad k\sum_{i=1}^{m} \lambda_{i,s}\lambda_{i,s+1}c_i(x_{s+1}) = k\sum_{i=1}^{m}\lambda_{i,s}/\lambda_{i,s+1}(\lambda_{i,s+1}c_i(x_{s+1}))^2.$$

Summing up (159) from $s = s_0$ to $s = N$ we have

$$d(\lambda^*) - d(\lambda_{s_0}) > d(\lambda_{N+1}) - d(\lambda_{s_0}) \geq k\sum_{s=s_0}^{N}\sum_{i=1}^{m}\lambda_{i,s}/\lambda_{i,s+1}(\lambda_{i,s+1}c_i(x_{s+1}))^2.$$

Keeping in mind (139) we obtain asymptotic complementarity condition

$$(160) \qquad \lim_{s\to\infty}(\lambda_s, c(x_s)) = 0.$$

Therefore

$$d(\lambda^*) = \lim_{s\to\infty} d(\lambda_s) = \lim_{s\to\infty}[f(x_s) - (\lambda_s, c(x_s))] = \lim_{s\to\infty} f(x_s),$$

that is

$$(161) \qquad \lim_{s\to\infty} f(x_s) = d(\lambda^*) = f(x^*).$$

From Slater condition follows boundedness of $L^*$. Therefore from concavity of $d$ follows boundedness $\Lambda(\lambda_0) = \{\lambda \in \mathbb{R}_+^m : d(\lambda) \geq d(\lambda_0)\}$. For any monotone increasing sequence $\{d(\lambda_s)\}_{s=s_0}^{\infty}$ follows boundedness $\Lambda_s = \{\lambda \in \mathbb{R}_+^m : d(\lambda) \geq d(\lambda_s)\}$ and $\Lambda_0 \supset ... \supset \Lambda_s \supset \Lambda_{s+1} \supset ... \supset L^*$. Therefore from (108) we have

$$(162) \qquad d_H(L^*, \partial\Lambda_s) > d_H(L^*, \partial\Lambda_{s+1}), \ s \geq s_0.$$

From (161) and (162) and continuity of $d$ follows

$$\lim_{s\to\infty} d_H(L^*, \partial\Lambda_s) = 0.$$

3) The dual sequence $\{\lambda_s\}_{s=0}^{\infty} \subset \Lambda(\lambda_0)$ is bounded, therefore there is a converging subsequence $\{\lambda_{s_l}\}_{l=1}^{\infty}$: $\lim_{l\to\infty}\lambda_{s_l} = \bar\lambda$.

Consider two subsets of indices $I_+ = \{i : \bar\lambda_i > 0\}$ and $I_0 = \{i : \bar\lambda_i = 0\}$. From the asymptotic complementarity (160) follows $\lim_{s\to\infty} c_i(x_s) = 0$, $i \in I_+$.

There exist such subsequence $\{s_l\}_{l=1}^{\infty}$ that for any $i \in I_0$ we have $\lim_{l\to\infty}\lambda_{i,s_l} = 0$, therefore without loosing the generality we can assume that

$$\lambda_{i,s_{l+1}} \leq 0.5\lambda_{i,s_l}, \ i \in I_0.$$

Using the update formula (153) we obtain

$$\lambda_{s_{l+1}}\prod_{s=s_l}^{s_{l+1}}(k\lambda_{i,s}c_i(x_s) + 1) = \lambda_{i,s_l} \geq 2\lambda_{i,s_{l+1}}, \ i \in I_0.$$

Invoking the arithmetic-geometric means inequality for $i \in I_0$ we obtain

$$\frac{1}{s_{l+1} - s_l}\sum_{s=s_l}^{s_{l+1}}(k\lambda_{i,s}c_i(x_s) + 1) \geq \left(\prod_{s=s_l}^{s_{l+1}}(k\lambda_{i,s}c_i(x_s) + 1)\right)^{\frac{1}{s_{l+1} - s_l}} \geq 2^{\frac{1}{s_{l+1} - s_l}}$$

or

$$\sum_{s=s_l}^{s_{l+1}}\lambda_{i,s}c_i(x_s) > 0, \ i \in I_0.$$

Using Jensen inequality and concavity $c_i$ we obtain

$$c_i(\bar{x}_{l+1}) = c_i\left(\sum_{s=s_l}^{s_{l+1}} \bar{\lambda}_{i,s} x_s\right) \geq \sum_{s=s_l}^{s_{l+1}} \bar{\lambda}_{i,s} c_i(x_s) > 0,$$

where $\bar{\lambda}_{i,s} = \lambda_{i,s}\left(\sum_{s=s_l}^{s_{l+1}} \lambda_{i,s}\right)^{-1} \geq 0$, $\sum_{s=s_l}^{s_{l+1}} \bar{\lambda}_{i,s} = 1$, $i \in I_0$. Keeping in mind $\lim_{s\to\infty} c_i(x_s) = 0$, $i \in I_+$ we conclude that the sequence $\{\bar{x}_{l+1}\}_{l=0}^{\infty}$ is asymptotically feasible, therefore it is bounded. Without loosing generality we can assume that $\lim_{l\to\infty} \bar{x}_l = \bar{x} \in \Omega$.

From convexity $f$ follows

$$f(\bar{x}_{l+1}) \leq \sum_{s=s_l}^{s_{l+1}} \bar{\lambda}_{i,s} f(x_s).$$

Therefore from (161) follows

$$f(\bar{x}) = \lim_{l\to\infty} f(\bar{x}_{l+1}) \leq \lim_{s\to\infty} f(x_s) = \lim_{s\to\infty} d(\lambda_s) = d(\lambda^*) = f(x^*).$$

Thus, $f(\bar{x}) = f(x^*)$, hence $d(\lambda^*) = d(\bar{\lambda})$ and $\bar{x} = x^*$, $\bar{\lambda} = \lambda^*$. □

The items 1) and 2 a) of Theorem 18 were proven by Matioli and Gonzaga (see Theorem 3.2 in [39]).

6.4. **Lagrangian Transformation and Affine Scaling method for LP.** Let $a \in \mathbb{R}^n, b \in \mathbb{R}^m$ and $A : \mathbb{R}^n \to \mathbb{R}^m$ are given. We consider the following LP problem

$$(163) \qquad x^* \in X^* = Argmin\{(a,x)|c(x) = Ax - b \geq 0\}$$

and the dual LP

$$(164) \qquad \lambda^* \in L^* = Argmin\{(b,\lambda)|r(\lambda) = A^T\lambda - a = 0, \lambda \in \mathbb{R}_+^m\}.$$

The LT $\mathcal{L} : \mathbb{R}^n \times \mathbb{R}^m \times \mathbb{R}_{++} \to \mathbb{R}$ for LP is defined as follows

$$(165) \qquad \mathcal{L}(x,\lambda,k) := (a,x) - k^{-1}\sum_{s=1}^{m} \psi(k\lambda_i c_i(x)),$$

where $c_i(x) = (Ax - b)_i = (a_i, x) - b_i$, $i = 1, \ldots, m$.

We assume that $X^* \neq \phi$ is bounded and so is the dual optimal set $L^*$.

The LT method generate primal - dual sequence $\{x_{s+1}, \lambda_{s+1}\}_{s=0}^{\infty}$ by the following formulas

$$(166) \qquad x_{s+1} : \nabla_x \mathcal{L}(x_{s+1}, \lambda_s, k_s) = 0$$

$$(167) \qquad \lambda_{s+1} : \lambda_{i,s+1} = \lambda_{i,s}\psi'(k_s\lambda_{i,s}c_i(x_{s+1})), i = 1, \ldots, m.$$

**Theorem 19.** *If the primal optimal $X^*$ is bounded, then the LT method (166)–(167) is well defined for any transformation $\psi \in \Psi$. For the dual sequence $\{\lambda_s\}_{s=0}^{\infty}$ generated by (167) the following statements hold true:*

1) *the LT method (166)–(167) is equivalent to the following Interior Prox*

   $$k(b, \lambda_{s+1}) - B_\varphi(\lambda_{s+1}, \lambda_s) = \max\{k(b, \lambda) - B_\varphi(\lambda, \lambda_s)|A^T\lambda = 0\},$$

   *where $B_\varphi(u,v) = \sum_{i=1}^{m} \varphi(\frac{u_i}{v_i})$ is the Bregman type distance;*

2) *there exists $s_0 > 0$ that for any $s \geq s_0$ the LT method with truncated MBF transformation $\psi_2(t)$ is equivalent to the affine scaling type method for the dual LP.*

**Proof**

1) We use the vector form for formula (167) assuming that the multiplication and division are componentwise, *i.e.* for vectors $a, b \in \mathbb{R}^n$, the vector $c = ab = (c_i = a_i b_i, \quad i = 1, \ldots, n)$ and the vector $d = a/b = (d_i = a_i/b_i, i = 1, \ldots, n)$. From (167) follows

$$(168) \qquad \frac{\lambda_{s+1}}{\lambda_s} = \psi'(k\lambda_s c(x_{s+1})).$$

Using again the inverse function formula we obtain

$$(169) \qquad k\lambda_s c(x_{s+1}) = \psi'^{-1}(\lambda_{s+1}/\lambda_s).$$

It follows from (166) and (167) that

$$\nabla_x \mathcal{L}(x_{s+1}, \lambda_s, k) = a - A^T \psi'(k\lambda_s c(x_{s+1}))\lambda_s = a - A^T \lambda_{s+1}$$
$$= \nabla_x L(x_{s+1}, \lambda_{s+1}) = 0,$$

therefore

$$d(\lambda_{s+1}) = L(x_{s+1}, \lambda_{s+1}) = (a, x_{s+1}) - (\lambda_{s+1}, Ax_{s+1} - b) \quad =$$
$$(a - A^T \lambda_{s+1}, x_{s+1}) + (b, \lambda_{s+1}) = (b, \lambda_{s+1}).$$

Using LEID $\psi'^{-1} = \psi^{*'}$ and $\varphi = -\psi^*$ we can rewrite (169) as follows

$$(170) \qquad -kc(x_{s+1}) - (\lambda_s)^{-1}\varphi'(\lambda_{s+1}/\lambda_s) = 0.$$

Keeping in mind $A^T \lambda_{s+1} = a, -c(x_{s+1}) \in \partial d(\lambda_{s+1})$ and $\lambda_{s+1} \in \mathbb{R}^m_{++}$ we can view (170) as the optimality criteria for the following problem

$$(171) \qquad k(b, \lambda_{s+1}) - B_\varphi(\lambda_{s+1}, \lambda_s) = \max\{kd(\lambda) - B_\varphi(\lambda, \lambda_s) | A^T \lambda = a\},$$

where $B_\varphi(\lambda, \lambda_s) = \sum_{i=1}^{q} \varphi(\lambda_i/\lambda_{i,s})$ is Bregman type distance.

2) Let's consider the LT method with truncated MBF transformation $\psi_2(t)$. It follows from (139) that there exists $s_0$ that for any $s \geq s_0$ only MBF kernel $\varphi_2 = -\ln s + s - 1$ and correspondent Bregman distance

$$\mathbb{B}_2(\lambda, \lambda_s) = \sum_{i=1}^{q}(-ln\frac{\lambda_i}{\lambda_{i,s}} + \frac{\lambda_i}{\lambda_{i,s}} - 1)$$

will be used in (171). Using considerations similar to those in item 4) Theorem 17 we can rewrite (171) as follows

$$(172) \qquad k(b, \lambda_{s+1}) = \arg\max\{k(b, \lambda) | \lambda \in E(\lambda_s, r_s), A^T \lambda = a\},$$

where $r_s^2 = Q(\lambda_{s+1}, \lambda_s)$ and $E(\lambda_s, r_s) = \{\lambda : (\lambda - \lambda_s)^T H_s^{-2}(\lambda - \lambda_s) \leq r_s^2\}$ is Dikin's ellipsoid and (172) is affine scaling type method for the dual LP (see [18]).

In the final part of the paper(survey) we will show the role of LET and LEINV in unconstrained minimization of SC functions. For the basic SC properties and damped Newton method see [42] and [43].

## 7. Legendre Invariant and Self-Concordant Functions

We consider a closed convex function $F \in C^3$ defined on an open convex set $\mathrm{dom}F \subset \mathbb{R}^n$. For a given $x \in \mathrm{dom}F$ and direction $u \in \mathbb{R}^n \setminus \{0\}$ we consider the restriction

$$f(t) = F(x + tu)$$

of $F$, which is defined on $\mathrm{dom}f = \{t : x + tu \in \mathrm{dom}F\}$. Along with $f$, let us consider its derivatives

$$\begin{aligned}
f'(t) &= (\nabla F(x + tu), u), \\
f''(t) &= (\nabla^2 F(x + tu)u, u), \\
f'''(t) &= (\nabla^3 F(x + tu)[u]u, u),
\end{aligned}$$

where $\nabla F$ is the gradient of $F$, $\nabla^2 F$ is the Hessian of $F$ and

$$\nabla^3 F(x)[u] = \lim_{\tau \to +0} \tau^{-1} \left[ \nabla^2 F(x + \tau u) - \nabla^2 F(x) \right].$$

Then,

$$\begin{aligned}
DF(x)[u] &:= (\nabla F(x), u) = f'(0), \\
D^2 F(x)[u, u] &:= (\nabla^2 F(x)u, u) = f''(0), \\
D^3 F(x)[u, u, u] &:= (\nabla^3 F(x)[u]u, u) = f'''(0).
\end{aligned}$$

Function $F$ is self-concordant if there is $M > 0$ such that the inequality

$$D^3 F(x)[u, u, u] \leq M(\nabla^2 F(x)u, u)^{\frac{3}{2}}$$

holds for any $x \in \mathrm{dom}F$ and any $u \in \mathbb{R}^n$.

If for a SC function $F$ the $\mathrm{dom}F$ does not contain a straight line, then the Hessian $\nabla^2 F(x)$ is positive definite at any $x \in \mathrm{dom}F$. We assume that such condition holds, so for any $x \in \mathrm{dom}F$ and any $u \in \mathbb{R}^n \setminus \{0\}$ we have

$$(173) \qquad (\nabla^2 F(x)u, u) = f''(0) > 0,$$

that is $F$ is strictly convex on $\mathrm{dom}F$.

A strictly convex function $F$ is self-concordant (SC) if the Legendre invariant of the restriction $f(t) = F(x + tu)$ is bounded, i.e. for any $x \in \mathrm{dom}F$ and any direction $u = y - x \in \mathbb{R}^n \setminus \{0\}$ there exist $M > 0$ that

$$(174) \qquad \mathrm{LEINV}(f) = |f'''(t)|(f''(t))^{-\frac{3}{2}} \leq M, \ \forall t : x + tu \in \mathrm{dom}F.$$

Let us consider the log-barrier function $F(x) = -\ln x$, then for any $x \in \mathrm{dom}F = \{x : x > 0\}$ we have $F'(x) = -x^{-1}$, $F''(x) = x^{-2}$, $F'''(x) = -2x^{-3}$ and

$$(175) \qquad \mathrm{LEINV}(F) = |F'''(x)| \ (F''(x))^{-3/2} \leq 2.$$

Therefore, $F(x) = -\ln x$ is self-concordant with $M = 2$.

The following function

$$g(t) = (\nabla^2 F(x + tu)u, u)^{-1/2} = (f''(t))^{-1/2},$$

is critical for the self-concordance (SC) theory.

For any $t \in \mathrm{dom}f$, we have

$$g'(t) = \frac{d \left[ (f''(t))^{-1/2} \right]}{dt} = -\frac{1}{2} f'''(t)(f''(t))^{-3/2}.$$

It follows from (175) that

$$(176) \qquad\qquad 0.5\, \mathrm{LEINV}(f) = |g'(t)| \leq 1\,, \quad \forall\, t \in \mathrm{dom}f.$$

The differential inequality (176) is the key element for establishing basic bounds for SC functions.

The other important component of the SC theory is two local scaled norms of a vector $u \in \mathbb{R}^n$. The first local scaled norm is defined at each point $x \in \mathrm{dom}F$ as follows

$$\|u\|_x = \left(\nabla^2 F(x)u, u\right)^{1/2}.$$

The second scaled norm is defined by formula

$$\|v\|_x^* = \left(\left(\nabla^2 F(x)\right)^{-1} v, v\right)^{1/2}.$$

From (173) follows that the second scaled norm is well defined at each $x \in \mathrm{dom}F$. The following Cauchy-Schwarz (CS) inequality for scaled norms will be often used later.

Let matrix $A \in \mathbb{R}^{n \times n}$ be symmetric and positive definite, then $A^{1/2}$ exists and

$$
\begin{aligned}
|(u,v)| &= \left|\left(A^{1/2}u,\, A^{-1/2}v\right)\right| \leq \left\|A^{1/2}u\right\| \left\|A^{-1/2}v\right\| \\
&= \left(A^{1/2}u,\, A^{1/2}u\right)^{1/2} \left(A^{-1/2}v,\, A^{-1/2}v\right)^{1/2} \\
&= (Au, u)^{1/2} \left(A^{-1}v, v\right)^{1/2} = \|u\|_A \ \|v\|_{A^{-1}}.
\end{aligned}
$$

By taking $A = \nabla^2 F(x)$, for any $u, v \in \mathbb{R}^n$ one obtains the following CS inequality:

$$|(u,v)| \leq \|u\|_x \ \|v\|_x^*.$$

The following Proposition will be used later.

**Proposition 20.** *A function $F$ is self-concordant if and only if for any $x \in \mathrm{dom}F$ and any $u_1, u_2, u_3 \in \mathbb{R}^n \setminus \{0\}$ we have*

$$(177) \qquad\qquad \left|D^3 F(x)\left[u_1, u_2, u_3\right]\right| \leq 2 \prod_{i=1}^{3} \|u_i\|_x\,,$$

*where $D^3 F(x)[u_1, u_2, u_3] = (\nabla^3 F(x)[u_1]u_2, u_3)$.*

The following theorem establishes one of the most important facts about SC functions: any SC function is a barrier function on $\mathrm{dom}F$. The opposite statement is, generally speaking, not true, that is not every barrier function is self-concordant. For example, the hyperbolic barrier $F(x) = x^{-1}$ defined on $\mathrm{dom}F = \{x : x > 0\}$ is not a SC function.

**Theorem 21.** *Let $F$ be a closed convex function on an open $\mathrm{dom}F$. Then, for any $\bar{x} \in \partial(\mathrm{dom}F)$ and any sequence $\{x_s\} \subset \mathrm{dom}F$ such that $x_s \to \bar{x}$, we have*

$$(178) \qquad\qquad \lim_{s \to \infty} F(x_s) = \infty\,.$$

*Proof.* From convexity of $F$ follows

$$F(x_s) \geq F(x_0) + (\nabla F(x_0),\, x_s - x_0)$$

for any given $x_0 \in \mathrm{dom}F$.

So, the sequence $\{F(x_s)\}$ is bounded from below. If (177) is not true, then the sequence $\{F(x_s)\}$ is bounded from above. Therefore, it has a limit point $\bar{F}$.

Without loss of generality, we can assume that $z_s = (x_s, F(x_s)) \to \bar{z} = (\bar{x}, \bar{F})$. Since the function $F$ is closed, we have $\bar{z} \in \mathrm{epi}F$, but it is impossible because $\bar{x} \notin \mathrm{dom}F$. Therefore for any sequence

$$\{x_s\} \subset \mathrm{dom}F : \lim_{s \to \infty} x_s = \bar{x} \in \partial(\mathrm{dom}F)$$

we have (177). It means that $F$ is a barrier function on the $cl(\mathrm{dom}F)$.   □

For any $x \in \mathrm{dom}F$, and any $u \in \mathbb{R}^n \setminus \{0\}$ from (173) follows

$$(179) \qquad \left(\nabla^2 F(x)u, u\right) = \|u\|_x^2 > 0$$

and for $\forall t \in \mathrm{dom}f$ we have

$$(180) \qquad g(t) = \left(\nabla^2 F(x + tu)u, u\right)^{-1/2} = \|u\|_{x+tu}^{-1} > 0.$$

7.1. **Basic Bounds for SC Functions.** In this section the basic bounds for SC functions will be obtained by integration of inequalities (176) and (177).

First Integration. Keeping in mind $f''(t) > 0$ from (176), for any $s > 0$, we obtain

$$-\int_0^s dt \leq \int_0^s d\left(f''(t)^{-1/2}\right) \leq \int_0^s dt \,.$$

Therefore

$$(181) \qquad f''(0)^{-1/2} - s \leq f''(s)^{-1/2} \leq f''(0)^{-1/2} + s$$

or

$$(182) \qquad \left(f''(0)^{-1/2} + s\right)^{-2} \leq f''(s) \leq \left(f''(0)^{-1/2} - s\right)^{-2} \,.$$

The left inequality in (182) holds for all $s \geq 0$, while the right inequality holds only for $0 \leq s < f''(0)^{-1/2}$.

Let $x, y \in \mathrm{dom}F$, $y \neq x$, $u = y - x$ and $y(s) = x + s(y - x)$, $0 \leq s \leq 1$, so $y(0) = x$ and $y(1) = y$. Therefore,

$$f''(0) = \left(\nabla^2 F(x)(y - x), \, y - x\right) = \|y - x\|_x^2$$

and

$$f''(0)^{1/2} = \|y - x\|_x \,.$$

Also,

$$f''(1) = \left(\nabla^2 F(y)(y - x), \, y - x\right) = \|y - x\|_y^2$$

and

$$f''(1)^{1/2} = \|y - x\|_y \,.$$

From (181), for $s = 1$ follows

$$f''(0)^{-1/2} - 1 \leq f''(1)^{-1/2} \leq f''(0)^{-1/2} + 1,$$

or

$$\frac{1}{\|y - x\|_x} - 1 \leq \frac{1}{\|y - x\|_y} \leq \frac{1}{\|y - x\|_x} + 1 \,.$$

From the right inequality follows

$$(183) \qquad \|y - x\|_y \geq \frac{\|y - x\|_x}{1 + \|y - x\|_x} \,.$$

If $\|y - x\|_x < 1$, then from the left inequality follows

$$(184) \qquad \|y - x\|_y \leq \frac{\|y - x\|_x}{1 - \|y - x\|_x} \,.$$

By integrating (176) we get

(185) $$g(t) \geqslant g(0) - |t| \, , \quad t \in \text{dom} f \, .$$

For $x + tu \in \text{dom} F$ from (180) follows $g(t) > 0$. From Theorem 21 follows $F(x + tu) \to \infty$ when $x + tu \to \partial(\text{dom} F)$. Therefore, $(\nabla^2 F(x+tu)u, u)$ cannot be bounded when $x + tu \to \partial(\text{dom} F)$. Therefore from (180) follows $g(t) \to 0$ when $x + tu \to \partial(\text{dom} F)$. It follows from (185) that any $t : |t| < g(0)$ belongs to $\text{dom} f$, i.e.

$$(-g(0), g(0)) = \left( -\|u\|_x^{-1}, \, \|u\|_x^{-1} \right) \subset \text{dom} f \, .$$

Therefore, the set

$$E^0(x,1) = \left\{ y = x + tu : t^2 \, \|u\|_x^2 < 1 \right\}$$

is contained in $\text{dom} F$. In other words, the *Dikin's ellipsoid*

$$E(x,r) = \left\{ y \in \mathbb{R}^n : \|y - x\|_x^2 \leq r \right\} \, ,$$

is contained in $\text{dom} F$ for any $x \in \text{dom} F$ and any $r < 1$.

One can expect that, for any $x \in \text{dom} F$ and any $y \in E(x,r)$, the Hessians $\nabla^2 F(x)$ and $\nabla^2 F(y)$ are "close" enough if $0 < r < 1$ is small enough. The second integration allows to establish the corresponding bounds.

Second Integration. Let us fix $x \in \text{dom} F$, for a given $y \in \text{dom} F$ $(y \neq x)$ consider direction $u = y - x \in \mathbb{R}^n \setminus \{0\}$. Let $y(t) = x + tu = x + t(y - x)$, then for $t \geq 0$ and $y(t) \in \text{dom} F$ we have

$$\psi(t) = \|u\|_{y(t)}^2 = (F''(y(t))u, u).$$

It follows from Proposition 20 that

$$|\psi'(t)| = D^3 F(y(t))[y - x, u, u] \leq 2 \, \|y - x\|_{y(t)} \, \|u\|_{y(t)}^2 = 2 \, \|y - x\|_{y(t)} \, \psi(t) \, .$$

First of all, $\|y(t) - x\|_x \leq \|y - x\|_x$ for any $t \in [0,1]$. Keeping in mind that $y - x = t^{-1}(y(t) - x)$ and assuming $\|y - x\|_x < 1$ from (184) follows

$$
\begin{aligned}
|\psi'(t)| &\leq \frac{2}{t} \, \|y(t) - x\|_{y(t)} \, \psi(t) \leq \frac{2}{t} \, \frac{\|y(t) - x\|_x}{1 - \|y(t) - x\|_x} \, \psi(t) \\
&\leq 2 \, \frac{\|y - x\|_x}{1 - t \, \|y - x\|_x} \, \psi(t) \, .
\end{aligned}
$$

Therefore for $0 < t < \|y - x\|_x^{-1}$ follows

$$\frac{|\psi'(t)|}{\psi(t)} \leq \frac{2 \, \|y - x\|_x}{1 - t \, \|y - x\|_x} \, .$$

By integrating the above inequality we get

$$-2 \int_0^s \frac{\|y - x\|_x}{1 - t \, \|y - x\|_x} \, dt \leq \int_0^s \frac{\psi'(t)}{\psi(t)} \, dt \leq 2 \int_0^s \frac{\|y - x\|_x}{1 - t \, \|y - x\|_x} \, dt \, ,$$

for $0 < s < \|y - x\|_x^{-1}$, hence

$$2 \ln \left( 1 - s \, \|y - x\|_x \right) \leq \ln \psi(s) - \ln \psi(0) \leq -2 \ln \left( 1 - s \, \|y - x\|_x \right) \, .$$

For $s = 1$, we have

(186) $$\psi(0) \left( 1 - \|y - x\|_x \right)^2 \leq \psi(1) \leq \psi(0) \left( 1 - \|y - x\|_x \right)^{-2} \, .$$

In view of $\psi(0) = (\nabla^2 F(x)u, u)$ and $\psi(1) = (\nabla^2 F(y)u, u)$ for any $u \in \mathbb{R}^n \setminus \{0\}$ from (186) follows

$$\left(1 - \|y - x\|_x\right)^2 \left(\nabla^2 F(x)u, u\right) \leq \left(\nabla^2 F(y)u, u\right) \leq \left(1 - \|y - x\|_x\right)^{-2} \left(\nabla^2 F(x)u, u\right).$$

Therefore the following matrix inequality holds

$$(187) \qquad \left(1 - \|y - x\|_x\right)^2 \nabla^2 F(x) \preccurlyeq \nabla^2 F(y) \preccurlyeq \nabla^2 F(x) \left(1 - \|y - x\|_x\right)^{-2},$$

where $A \succcurlyeq B$ means that $A - B$ is nonnegative definite. Note that (187) takes place for any $x, y \in \mathrm{dom} F$.

In order to find the upper and the lower bounds for the matrix

$$(188) \qquad G = \int_0^1 \nabla^2 F(x + \tau(y - x)) d\tau$$

let us consider (187) for $y := x + \tau(y - x)$.

From the left inequality (187) follows

$$G = \int_0^1 \nabla^2 F(x + \tau(y - x)) d\tau \succcurlyeq \nabla^2 F(x) \int_0^1 \left(1 - \tau \|y - x\|_x\right)^2 d\tau.$$

Therefore, for $r = \|y - x\|_x < 1$, we have

$$(189) \qquad G \succcurlyeq \nabla^2 F(x) \int_0^1 (1 - \tau r)^2 d\tau = \nabla^2 F(x) \left(1 - r + \frac{r^2}{3}\right).$$

From the right inequality (187) follows

$$(190) \qquad G \preccurlyeq \nabla^2 F(x) \int_0^1 (1 - \tau r)^{-2} d\tau = \nabla^2 F(x) \frac{1}{1 - r},$$

i.e. for any $x \in \mathrm{dom} F$, the following inequalities hold:

$$(191) \qquad \left(1 - r + \frac{r^2}{3}\right) \nabla^2 F(x) \preccurlyeq G \preccurlyeq \frac{1}{1 - r} \nabla^2 F(x).$$

The first two integrations produced two very important facts.

(1) For any $x \in \mathrm{dom} F$, Dikin's ellipsoid

$$E(x, r) = \left\{y \in \mathbb{R}^n : \|y - x\|_x^2 \leq r\right\}$$

is contained in $\mathrm{dom} F$, for any $0 \leq r < 1$.

(2) For any $x \in \mathrm{dom} F$ and any $y \in E(x, r)$ from (187) follows

$$(192) \qquad (1 - r)^2 \nabla^2 F(x) \preccurlyeq \nabla^2 F(y) \preccurlyeq \frac{1}{(1 - r)^2} \nabla^2 F(x),$$

i.e. the function $F$ is almost quadratic inside the ellipsoid $E(x, r)$ for small $0 \leq r < 1$.

The bounds for the gradient $\nabla F(x)$, which is a monotone operator in $\mathbb{R}^n$, we establish by integrating (182).

**Third Integration.** From (182), for $0 \leq t < f(0)^{-1/2} = \|y-x\|_x^{-1}$ and $0 \leq s \leq 1$ we obtain

$$\int_0^s \left( f''(0)^{-1/2} + t \right)^{-2} dt \leq \int_0^s f''(t) dt \leq \int_0^s \left( f''(0)^{-1/2} - t \right)^{-2} dt \,,$$

or

$$\begin{aligned}
f'(0) \quad + \quad & f''(0)^{1/2} \left( 1 - \left( 1 + s f''(0)^{1/2} \right)^{-1} \right) \\
\leq \quad & f'(s) \leq f'(0) - f''(0)^{1/2} \left( 1 - \left( 1 - s f''(0)^{1/2} \right)^{-1} \right) \,.
\end{aligned}$$

The obtained inequalities we can rewrite as follows

$$(193) \qquad f'(0) + w'(f''(0)^{\frac{1}{2}} s) \leq f'(s) \leq f'(0) + w^{*'}(f''(0)^{\frac{1}{2}} s),$$

where $\omega(t) = t - \ln(1+t)$ and $\omega^*(s) = \sup_{t>-1}\{st - t + \ln(1+t)\} = -s - \ln(1-s) = \omega(-s)$ is the LET of $\omega(t)$.

From the right inequality (193), for $s = 1$ follows

$$f'(1) - f'(0) \leq -f''(0)^{1/2} \left( 1 - \frac{1}{1 - f''(0)^{1/2}} \right) = \frac{f''(0)}{1 - f''(0)^{1/2}} \,.$$

Recalling formulas for $f'(0)$, $f'(1)$, $f''(0)$, and $f''(1)$ we get

$$(194) \qquad (\nabla F(y) - \nabla F(x), y - x) \leq \frac{\|y-x\|_x^2}{1 - \|y-x\|_x}$$

for any $x$ and $y \in \mathrm{dom} F$.

From the left inequality in (193), for $s = 1$ follows

$$f'(1) - f'(0) \geq f''(0)^{1/2} \left( 1 - \frac{1}{1 + f''(0)^{1/2}} \right) = \frac{f''(0)}{1 + f''(0)^{1/2}}$$

or

$$(195) \qquad (\nabla F(y) - \nabla F(x), y - x) \geq \frac{\|y-x\|_x^2}{1 + \|y-x\|_x} \,.$$

**Fourth Integration.** In order to establish bounds for $F(y) - F(x)$ it is enough to integrate the inequalities (193).

Taking the integral of the right inequality (193), we obtain

$$\begin{aligned}
f(s) \quad \leq \quad & f(0) + f'(0)s + \omega^* \left( f''(0)^{1/2} s \right) \\
= \quad & f(0) + f'(0)s - f''(0)^{1/2} s - \ln \left( 1 - f''(0)^{1/2} s \right) \\
(196) \qquad = \quad & U(s) \,.
\end{aligned}$$

In other words, $U(s)$ is an upper bound for $f(s)$ on the interval $[0, f''(0)^{-1/2})$. Recall that $f''(0)^{-1/2} = \|y-x\|_x^{-1} > 1$. For $s = 1$ from (196) follows

$$(197) \qquad f(1) - f(0) \leq f'(0) + \omega^* \left( f''(0)^{1/2} \right) = f'(0) + \omega^* \left( \|y-x\|_x \right) \,.$$

Keeping in mind $f(0) = F(x)$, $f(1) = F(y)$, from (197), we get

$$(198) \qquad F(y) - F(x) \leq (\nabla F(x), y - x) + \omega^* \left( \|y-x\|_x \right) \,.$$

Integration of the left inequality (193) leads to the lower bound $L(s)$ for $f(s)$

$$
\begin{aligned}
f(s) &\geq f(0) + f'(0)s + \omega\left(f''(0)^{1/2}s\right) \\
&= f(0) + f'(0)s + f''(0)^{1/2}s - \ln\left(1 + f''(0)^{1/2}s\right) \\
&= L(s)\,, \ \forall s \geq 0\,.
\end{aligned}
$$

(199)

For $s = 1$, we have

$$
f(1) - f(0) \geqslant f'(0) + \omega\left(f''(0)^{1/2}\right)
$$

or

(200)
$$
F(y) - F(x) \geq (\nabla F(x), y - x) + \omega\left(\|y - x\|_x\right)\,.
$$

We conclude the section by considering the existence of the minimizer

(201)
$$
x^* = \arg\min\{F(x) : x \in \operatorname{dom} F\}
$$

for a self-concordant function $F$.

It follows from (173) that the Hessian $\nabla^2 F(x)$ is positive definite for any $x \in \operatorname{dom} F$, but the existence of $x^* : \nabla F(x^*) = 0$, does not follow from strict convexity of $F$.

However, it guarantees the existence of the local norm $\|v\|_x^* = \left(\left(\nabla^2 F(x)\right)^{-1} v, v\right)^{1/2} > 0$ at any $x \in \operatorname{dom} F$.

For $v = \nabla F(x)$, one obtains the following scaled norm of the gradient $\nabla F(x)$,

$$
\lambda(x) = \left(\nabla^2 F(x)^{-1} \nabla F(x),\, \nabla F(x)\right)^{1/2} = \|\nabla F(x)\|_x^* > 0\,,
$$

which plays an important role in SC theory. It is called Newton decrement of $F$ at the point $x \in \operatorname{dom} F$.

**Theorem 22.** *If $\lambda(x) < 1$ for some $x \in \operatorname{dom} F$ then the minimizer $x^*$ in (201) exists.*

*Proof.* For $u = y - x \neq 0$ and $v = \nabla F(x)$, where $x$ and $y \in \operatorname{dom} F$ from CS inequality $|(u, v)| \leq \|v\|_x^* \|u\|_x$ follows

(202)
$$
|(\nabla F(x), y - x)| \leq \|\nabla F(x)\|_x^* \|y - x\|_x\,.
$$

From (200) and (202) and the formula for $\lambda(x)$ follows

$$
F(y) - F(x) \geq -\lambda(x) \|y - x\|_x + \omega\left(\|y - x\|_x\right)\,.
$$

Therefore, for any $y \in \mathcal{L}(x) = \{y \in \mathbb{R}^n : F(y) \leq F(x)\}$ we have

$$
\omega\left(\|y - x\|_x\right) \leq \lambda(x) \|y - x\|_x\,,
$$

i.e.

$$
\|y - x\|_x^{-1} \omega\left(\|y - x\|_x\right) \leq \lambda(x) < 1\,.
$$

From the definition of $\omega(t)$ follows

$$
1 - \frac{1}{\|y - x\|_x} \ln\left(1 + \|y - x\|_x\right) \leq \lambda(x) < 1\,.
$$

The function $1 - \tau^{-1}\ln(1 + \tau)$ is monotone increasing for $\tau > 0$. Therefore, for a given $0 < \lambda(x) < 1$, the equation

$$
1 - \lambda(x) = \tau^{-1}\ln(1 + \tau)
$$

has a unique root $\bar{\tau} > 0$. Thus, for any $y \in \mathcal{L}(x)$, we have

$$\|y - x\|_x \leq \bar{\tau},$$

i.e. the level set $\mathcal{L}(x)$ at $x \in \text{dom}F$ is bounded and closed due to the continuity of $F$. Therefore, $x^*$ exists due to the Weierstrass theorem. The minimizer $x^*$ is unique due to the strict convexity of $F(x)$ for $x \in \text{dom } F$. □

The theorem presents an interesting result: a local condition $\lambda(x) < 1$ at some $x \in \text{dom}F$ guarantees the existence of $x^*$, which is a global property of $F$ on the $\text{dom}F$. The condition $0 < \lambda(x) < 1$ will plays an important role later.

Let us briefly summarize the basic properties of the SC functions established so far.

(1) The SC function $F$ is a barrier function on $\text{dom}F$.
(2) For any $x \in \text{dom}F$ and any $0 < r < 1$, there is a Dikin's ellipsoid inside $\text{dom}F$, i.e.

$$E(x, r) = \left\{ y : \|y - x\|_x^2 \leq r \right\} \subset \text{dom}F.$$

(3) For any $x \in \text{dom}F$ and small enough $0 < r < 1$, the function $F$ is almost quadratic inside of the Dikin's ellipsoid $E(x, r)$ due to the bounds (192).
(4) The gradient $\nabla F$ is a strictly monotone operator on $\text{dom}F$ with upper and lower monotonicity bounds given by (194) and (195).
(5) For any $x \in \text{dom}F$ and any direction $u = y - x$, the restriction $f(s) = F(x + s(y - x))$ is bounded by $U(s)$ and $L(s)$ (see (196) and (199)).
(6) Condition $0 < \lambda(x) < 1$ at any $x \in \text{dom}F$ guarantees the existence of a unique minimizer $x^*$ on $\text{dom}F$.

It is quite remarkable that practically all important properties of SC functions follow from a single differential inequality (176), which is, a direct consequence of the boundedness of LEINV($f$).

We conclude the section by showing that Newton method can be very efficient for global minimization of SC functions, in spite of the fact that $F$ is not strongly convex.

### 7.2. Damped Newton Method for Minimization of SC Function.

The SC functions are strictly convex on $\text{dom}F$. Such a property, generally speaking, does not guarantee global convergence of the Newton method. For example, $f(t) = \sqrt{1 + t^2}$ is strictly convex, but Newton method for finding $\min_t f(t)$ diverges from any starting point $t_0 \notin ]-1, 1[$.

Turns out that SC properties guarantee convergence of the special damped Newton method from any starting point. Moreover, such method goes through three phases. In the first phase each step reduces the error bound $\Delta f(x) = f(x) - f(x^*)$ by a constant, which is independent on $x \in \text{dom}F$. In the second phase the error bound converges to zero with at least superlinear rate. The superlinear rate is characterized explicitly through $w(\lambda)$ and its LET $w^*(\lambda)$, where $0 < \lambda < 1$ is the Newton decrement. At the final phase the damped Newton method practically turns into standard Newton method and the error bound converges to zero with quadratic rate.

The following bounds for the restriction $f(s) = F(x + su)$ at $x \in \text{dom}F$ in the direction $u = y - x \in \mathbb{R}^n \setminus \{0\}$ is our main tool

$$(203) \qquad \mathcal{L}(s) \leq f(s) \leq U(s) \atop s \geq 0 \qquad\qquad 0 \leq s \leq f''(0)^{-(1/2)}.$$

Let $x \in \operatorname{dom} F$, $f(0) = F(x)$ and $x \neq x^*$, then there exists $y \in \operatorname{dom} F$ such that for $u = y - x \neq 0$ we have

$$\begin{aligned} a) \qquad & f'(0) = (\nabla F(x), u) < 0, \text{ and} \\ (204) \qquad b) \qquad & f''(0) = (\nabla^2 F(x)u, u) = \|u\|_x^2 = d^2 > 0. \end{aligned}$$

We would like to estimate the reduction of $F$, as a result of one Newton step with $x \in \operatorname{dom} F$ as a starting point.

Let us consider the upper bound

$$U(s) = f(0) + f'(0)s - ds - \ln(1 - ds),$$

for $f(s)$. The function $U(s)$ is strongly convex in $s$ on $[0, d)^{-1}$. Also, $U'(0) = f'(0) < 0$ and $U'(s) \to \infty$ for $s \to d^{-1}$. Therefore, the equation

$$(205) \qquad U'(s) = f'(0) - d + d(1 - ds)^{-1} = 0$$

has a unique solution $\bar{s} \in [0, d^{-1})$, which is the unconstrained minimizer for $U(s)$. From (205) we have

$$\bar{s} = -f'(0)d^{-2} \left(1 - f'(0)d^{-1}\right)^{-1} = \Delta(1 + \lambda)^{-1}$$

where $\Delta = -f'(0)d^{-2}$ and $0 < \lambda = -f'(0)d^{-1} < 1$. On the other hand, the unconstrained minimizer $\bar{s}$ is a result of one step of the damped Newton method for finding $\min_{s \geq 0} U(s)$ with step length $t = (1 + \lambda)^{-1}$ from $s = 0$ as a starting point. It is easy to see that

$$U\left((1 + \lambda)^{-1}\Delta\right) = f(0) - \omega(\lambda).$$

From the right inequality in (203), we obtain

$$(206) \qquad f\left((1 + \lambda)^{-1}\Delta\right) \leq f(0) - \omega(\lambda).$$

Keeping in mind (204) for the Newton direction $u = y - x = -(\nabla^2 F(x))^{-1}\nabla F(x)$ we obtain

$$\Delta = -\frac{f'(0)}{f''(0)} = -\frac{(\nabla F(x), u)}{(\nabla^2 F(x)u, u)} = 1.$$

In view of $f(0) = F(x)$, we can rewrite (206) as follows:

$$(207) \qquad F\left(x - (1 + \lambda)^{-1}(\nabla^2 F(x))^{-1}\nabla F(x)\right) \leq F(x) - \omega(\lambda).$$

In other words, finding an unconstrained minimizer of the upper bound $U(s)$ is equivalent to one step of the damped Newton method

$$(208) \qquad x_{k+1} = x_k - (1 + \lambda(x_k))^{-1}\left(\nabla^2 F(x_k)\right)^{-1}\nabla F(x_k)$$

for minimization of $F(x)$ on $\operatorname{dom} F$. Moreover, our considerations are independent from the starting point $x \in \operatorname{dom} F$. Therefore, for any starting point $x_0 \in \operatorname{dom} F$ and $k \geq 1$, we have

$$(209) \qquad F(x_{k+1}) \leq F(x_k) - \omega(\lambda).$$

The bound (209) is universal, i.e. it is true for any $x_k \in \operatorname{dom} F$.

Let us compute $\lambda = f'(0)f''(0)^{-1/2}$ for the Newton direction

$$u = -\nabla^2 F(x)^{-1}\nabla F(x).$$

We have

$$
\begin{aligned}
\lambda \equiv \lambda(x) &= -f'(0)\, f''(0)^{-1/2} \\
&= -\frac{(\nabla F(x), u)}{(\nabla^2 F(x)u, u)^{1/2}} \\
&= \left(\nabla^2 F(x)^{-1}\nabla F(x), \nabla F(x)\right)^{1/2} \\
&= \|\nabla F(x)\|_x^* \,.
\end{aligned}
$$

We have seen already that it is critical that $0 < \lambda(x_k) < 1, \ \ \forall k \geq 0$.

The function $\omega(t) = t - \ln(1+t)$ is a monotone increasing, therefore for a small $\beta > 0$ and $1 > \lambda(x) \geq \beta$, from (209) we obtain reduction of $F(x)$ by a constant $\omega(\beta)$ at each damped Newton step. Therefore, the number of damped Newton steps is bounded by

$$
N \leq (\omega(\beta))^{-1}(F(x^0) - F(x^*)) \,.
$$

The bound (209), however, can be substantially improved for

$$
x \in S(x^*, r) = \{x \in \operatorname{dom} F : F(x) - F(x^*) \leq r\}
$$

and $0 < r < 1$.

Let us consider the lower bound

$$
L(s) = f(0) + f'(0)s + ds - \ln(1 + ds) \leq f(s), \quad s \geq 0 \,.
$$

The function $L(s)$ is strictly convex on $s \geq 0$. If $0 < \lambda = -f'(0)d^{-1} < 1$, then

$$
L'\left(\Delta(1-\lambda)^{-1}\right) = 0 \,.
$$

Therefore,

$$
\bar{\bar{s}} = \Delta(1-\lambda)^{-1} = \arg\min\{L(s) \mid s \geq 0\}
$$

and

$$
L(\bar{\bar{s}}) = f(0) - \omega(-\lambda) \,.
$$

Along with $\bar{s}$ and $\bar{\bar{s}}$ we consider (see Fig. 2)

$$
s^* = \operatorname{argmin}\{f(s) \mid s \geq 0\} \,.
$$

For a small $0 < r < 1$ and $x \in S(x^*, r)$, we have $f(0) - f(s^*) < 1$, hence $f(0) - f(\bar{s}) < 1$. The relative progress per step is more convenient to measure on the logarithmic scale

$$
\kappa = \frac{\ln(f(\bar{s}) - f(s^*))}{\ln(f(0) - f(s^*))} \,.
$$

From $\omega(\lambda) < f(0) - f(s^*) < 1$ follows $-\ln \omega(\lambda) > -\ln(f(0) - f(s^*))$ or $\ln(f(0) - f(s^*)) > \ln \omega(\lambda)$. From $f(\bar{s}) \leq f(0) - \omega(\lambda)$ and $f(s^*) \geq f(0) - \omega(-\lambda)$ follows (see Fig 2)

$$
f(\bar{s}) - f(s^*) \leq \omega(-\lambda) - \omega(\lambda) \,.
$$

Hence,

$$
\ln(f(\bar{s}) - f(s^*)) < \ln(\omega(-\lambda) - \omega(\lambda))
$$

and

$$
\begin{aligned}
\kappa(\lambda) &\leq \frac{\ln(\omega(-\lambda) - \omega(\lambda))}{\ln \omega(\lambda)} \\
&= \frac{\ln\left(-2\lambda + \ln(1+\lambda)(1-\lambda)^{-1}\right)}{\ln(\lambda - \ln(1+\lambda))} \,.
\end{aligned}
$$

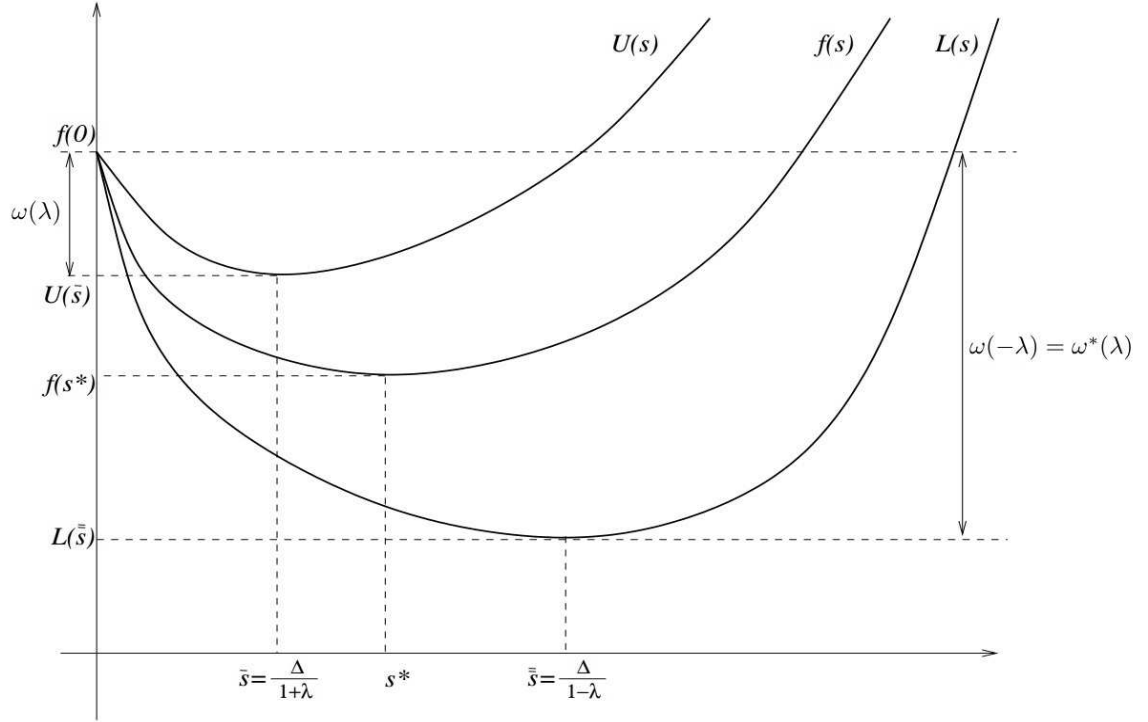

FIGURE 2.

For $0 < \lambda \le 0.5$, we have

$$\kappa(\lambda) \le \frac{\ln\left(\frac{2\lambda^3}{3} + \frac{2\lambda^5}{5}\right)}{\ln\left(\frac{\lambda^2}{2} - \frac{\lambda^3}{3} + \frac{\lambda^4}{4}\right)}.$$

In particular, $\kappa(0.40) \approx 1.09$. Thus, the sequence $\{x_k\}_{k=0}^{\infty}$ generated by the damped Newton method (208) with $\lambda(x_k) = 0.40$ converges in value at least with 1.09 $Q$-superlinear rate, that is for the error bound the $\Delta(x_k) = F(x_k) - F(x^*) < 1$, we have $\Delta(x_{k+1}) \le (\Delta(x_k))^{1.09}$.

Due to $\lim_{k\to\infty} \lambda(x_k) = \lim_{k\to\infty} \|\nabla F(x_k)\|_x = 0$ from some point on, method (208) practically turns into the classical Newton method

$$(210) \qquad x_{k+1} = x_k - \nabla^2 F(x_k)^{-1} \nabla F(x_k),$$

which converges with quadratic rate.

Instead of waiting for this to happen, there is a way of switching, at some point, from (208) to (210) and guarantee that from this point on, only Newton method (210) is used. Using such a strategy it is possible to achieve quadratic convergence earlier.

The following Theorem characterize the neighborhood at $x^*$ when quadratic convergence accuracy.

**Theorem 23.** *Let $x \in \mathrm{dom}F$ and*

$$\lambda(x) = \left(\nabla^2 F(x)^{-1}\nabla F(x), \nabla F(x)\right)^{1/2} < 1\,,$$

*then,*

    (1) *the point*

(211) $$\hat{x} = x - \nabla^2 F(x)^{-1}\nabla F(x)$$

        *belongs to* $\mathrm{dom}F$;

    (2) *the following bound holds*

(212) $$\lambda(\hat{x}) \leq \left(\frac{\lambda(x)}{1-\lambda(x)}\right)^2\,.$$

*Proof.* 1. Let $p = \hat{x} - x = -\nabla^2 F(x)^{-1}\nabla F(x)$, $\lambda = \lambda(x)$, then

$$
\begin{aligned}
\|p\|_x &= \left(\nabla^2 F(x)p, p\right)^{1/2} = \left(\nabla F(x),\, \nabla^2 F(x)^{-1}\nabla F(x)\right)^{1/2}\\
&= \|\nabla F(x)\|_x^* = \lambda(x) = \lambda < 1\,;
\end{aligned}
$$

therefore, $\hat{x} \in \mathrm{dom}F$.

2. First of all, note that if $A = A^T \succ 0$, $B = B^T \succ 0$ and $A \succcurlyeq B$, then

$$A^{-1} - B^{-1} = -A^{-1}(A-B)B^{-1} \preccurlyeq 0\,.$$

For $y = \hat{x}$ from the left inequality in (187), we obtain

$$
\begin{aligned}
\lambda(\hat{x}) &= \|\nabla F(\hat{x})\|_{\hat{x}}^* \leq (1 - \|p\|_x)^{-1}\left(\nabla^2 F(x)^{-1}\nabla F(\hat{x}),\, \nabla F(\hat{x})\right)^{1/2}\\
&= (1 - \|p\|_x)^{-1}\|\nabla F(\hat{x})\|_x^*\,.
\end{aligned}
$$

We can then rewrite (211) as follows

$$\nabla^2 F(x)\,(\hat{x} - x) + \nabla F(x) = 0\,.$$

Therefore,

$$\nabla F(\hat{x}) = \nabla F(\hat{x}) - \nabla F(x) - \nabla^2 F(x)(\hat{x} - x)\,.$$

Then, using (188) and formula (13) (see p. 6 [ ]), we obtain

$$\nabla F(\hat{x}) - \nabla F(x) = \int_0^1 \left(\nabla^2 F(x + \tau(\hat{x} - x))\,(\hat{x} - x)d\tau = G(\hat{x} - x)\,.\right.$$

Hence,

$$\nabla F(\hat{x}) = \left(G - \nabla^2 F(x)\right)(\hat{x} - x) = \hat{G}(\hat{x} - x) = \hat{G}p$$

and $\hat{G}^T = \hat{G}$.

From CS inequality follows

$$
\begin{aligned}
\|\nabla F(\hat{x})\|_x^{*2} &= \left(\nabla^2 F(x)^{-1}\hat{G}p, \hat{G}p\right) = \left(\hat{G}\nabla^2 F(x)^{-1}\hat{G}p, p\right)\\
(213) &\leq \left\|\hat{G}\nabla^2 F(x)^{-1}\hat{G}p\right\|_x^* \|p\|_x\,.
\end{aligned}
$$

Then

$$
\begin{aligned}
\left\|\hat{G}\nabla^2 F(x)^{-1}\hat{G}p\right\|_x^* &= \left(\hat{G}\nabla^2 F(x)^{-1}\hat{G}p,\ \nabla^2 F(x)^{-1}\hat{G}\nabla^2 F(x)^{-1}\hat{G}p\right)^{1/2} \\
&= \left(H(x)^2\nabla^2 F(x)^{-1/2}\hat{G}p,\ \nabla^2 F(x)^{-1/2}\hat{G}p\right)^{1/2} \\
&\leq \|H(x)\|\left(\nabla^2 F(x)^{-1/2}\hat{G}p, \nabla^2 F(x)^{-1/2}\hat{G}p\right)^{1/2} \\
&= \|H(x)\|\left(\nabla^2 F(x)^{-1}\hat{G}p, \hat{G}p\right) \\
&= \|H(x)\|\left(\nabla^2 F(x)^{-1}\nabla F(\hat{x}), \nabla F(\hat{x})\right)^{1/2} \\
&= \|H(x)\|\ \|\nabla F(\hat{x})\|_x^*\,,
\end{aligned}
$$

where $H(x) = \nabla^2 F(x)^{-1/2}\hat{G}\nabla^2 F(x)^{-1/2}$, therefore $\nabla^2 F(x)^{\frac{1}{2}}H(x)\nabla^2 F^{\frac{1}{2}}(x) = \hat{G}$.

From (213) and the last inequality we obtain

$$
\|\nabla F(\hat{x})\|_x^* \leq \|H(x)\|\ \|p\|_x = \lambda\|H(x)\|\,.
$$

It follows from (191)

$$
\left(-\lambda + \frac{\lambda^2}{3}\right)\nabla^2 F(x) \preccurlyeq \hat{G} = G - \nabla^2 F(x) \preccurlyeq \frac{\lambda}{1-\lambda}\,\nabla^2 F(x)\,.
$$

Then,

$$
\|H(x)\| \leq \max\left\{\frac{\lambda}{1-\lambda}\,,\, -\lambda+\frac{\lambda^2}{3}\right\} = \frac{\lambda}{1-\lambda}\,.
$$

Therefore,

$$
\lambda^2(\hat{x}) \leq \frac{1}{(1-\lambda)^2}\,\|\nabla F(\hat{x})\|_x^{*2} \leq \frac{1}{(1-\lambda)^2}\,\lambda^2\|H(x)\|^2 \leq \frac{\lambda^4}{(1-\lambda)^4}
$$

or

$$
\lambda(\hat{x}) \leq \frac{\lambda^2}{(1-\lambda)^2}\,.
$$

We saw already that $\lambda = \lambda(x) < 1$ is the main ingredient for the damped Newton method (208) to converge. To retain the same condition for $\lambda(\hat{x})$, it is sufficient to require $\lambda(\hat{x}) \leq \lambda \leq \lambda^2/(1-\lambda)^2$. The function $[\lambda/(1-\lambda)]^2$ is positive and monotone increasing on $(0,1)$. Therefore, to find an upper bound for $\lambda$ it is enough to solve the equation $\lambda/(1-\lambda)^2 = 1$. In other words, for any $\lambda = \lambda(x) < \bar{\lambda} = \frac{3-\sqrt{5}}{2}$, we have

$$
\lambda(\hat{x}) \leq \left(\frac{\lambda}{1-\lambda}\right)^2\,.
$$

Thus, the damped Newton method (208) follows three major stages in terms of the rate of convergence. First, it reduces the function value by a constant at each step. Then, it converges with superlinear rate and ,finally, in the neighborhood of the solution it converges with quadratic rate.

The Newton area, where the Newton method converges with the quadratic rate is defined as follows:

$$
(214)\qquad N(x^*, \beta) = \left\{x : \lambda(x) = \|\nabla F(x)\|_x^* \leq \beta < \bar{\lambda} = \frac{3-\sqrt{5}}{2}\right\}\,.
$$

To speed up the damped Newton method (208) one can use the following switching strategy. For a given $0 < \beta < \bar{\lambda} = (3 - \sqrt{5})/2$, one uses the damped Newton method (208) if $\lambda(x_k) > \beta$ and the "pure" Newton method (210) when $\lambda(x_k) \leq \beta$.

## 8. Concluding Remarks

The LEID is an universal instrument for establishing the duality results for SUMT, NR and LT methods. The duality result, in turn, are critical for both understanding the convergence mechanisms and the convergence analysis.

In particular, the update formula (107) and concavity of the dual function $d$ leads to the following bound

$$d(\lambda_{s+1}) - d(\lambda_s) \geq (kL)^{-1} \|\lambda_{s+1} - \lambda_s\|^2,$$

which together with $d(\lambda_{s+1}) - d(\lambda_s) \to 0$ shows that the Lagrange multipliers do not change much from same point on. It means that if Newton method is used for primal minimization then, from some point on, usually after very few Lagrange multipliers update the approximation for the primal minimizer $x_s$ is in the Newton area for the next minimizer $x_{s+1}$.

Therefore it takes few and, from some point on, only one Newton step to find the next primal approximation and update the Lagrange multipliers.

This phenomenon is called - the "hot" start (see [46]). The neighborhood of the solution where the "hot" start occurs has been characterized in [37] and observed in [5], [10], [25], [41].

It follows from Remark 14 that, under standard second order optimality condition, each Lagrange multipliers update shrinks the distance between the current and the optimal solution by a factor, which can be made as small as one wants by increasing $k > 0$.

In contrast to SUMT the NR methods requires much less computational effort per digit of accuracy at the end of the process then at the beginning.

Therefore NR methods is used when high accuracy needed (see, for example, [1]).

One of the most important features of NR methods is their numerical stability. It is due to the stability of the Newton's area, which does not shrink to a point in the final phase. Therefore one of the most reliable NLP solver PENNON is based on NR methods (see [32]-[34]).

The NR method with truncated MBF transformation has been widely used for both testing the NLP software and solving real life problems (see [1], [5], [10], [25], [32]-[34], [37] , [41] ). The numerical results obtained strongly support the theory, including the "hot" start phenomenon.

The NR as well as LT are primal exterior points methods. Their dual equivalence are interior points methods.

In particular, the LT with MBF transform $\psi(t) = \ln(t + 1)$ leads to the interior prox with Bregman distance, which is based on the self-concordant MBF kernel $\varphi(s) = -\psi^*(s) = -\ln s + s - 1$. Application of this LT for LP calculations leads to Dikin's type interior point method for the dual LP. It establishes, eventually, the remarkable connection between exterior and interior point methods (see [39], [49]).

On the other hand, the LEINV is in the heart of the SC theory - one of the most beautiful chapters of the modern optimization.

Although the Legendre Transformation was introduced more than 200 years ago, we saw that LEID and LEINV are still critical in modern optimization both constrained and unconstrained.

## References

[1] Alber M., Reemtsen R.: Intensity modulated radiotherapy treatment planning by use of a barrier-penalty multiplier method. Optimization Methods and Software. 22, N 3, 391-411, (2007)

[2] Antipin A. S.: Methods of nonlinear programming based on the direct and dual augmentation of the Lagrangian. Moscow VNIISI, (1979)

[3] Auslender R., Cominetti R. and Haddou M.: Asymptotic analysis for penalty and barrier methods in convex and linear programming, Mathematics of Operations Research 22(1), 43-62, (1997)

[4] Bauschke H., Matouskova E. and Reich S.: Projection and proximal point methods, convergence results and counterexamples. Nonlinear Anal. 56, no. 5, 715-738, (2004)

[5] Ben-Tal A., Nemirovski A.: Optimal design of Engineering Structures. Optima, 47, 4-9 (1995)

[6] Ben-Tal A., Zibulevski M.: Penalty-barrier methods for convex programming problems SIAM J. Optim. 7, 347-366 (1997)

[7] Bertsekas D.: Constrained optimization and Lagrange multiplier methods. New York, (1982)

[8] Bregman L.: The relaxation method for finding the common point of convex sets and its application to the solution of problems in convex programming," USSR Computational Mathematics and Mathematical Physics, 7, 200–217, (1967)

[9] Bregman L., Censor Y., Reich S.: Dykstra Algorithm as the Nonlinear Extention of Bregman's optimization Method. Journal of Convex Analysis, 6, N2, 319-333, (1999)

[10] Breitfeld M., Shanno D.: "Computational experience with modified log-barrier methods for nonlinear programming, Annals of OR, 62, 439-464 (1996)

[11] Byrne C., Censor Y.: Proximity Function minimization Using Multiple Bregman Projections with Application to Split Feasibility and Kullback-Leibler Distance Minimization. Annals of OR, 105, 77-98, (2001)

[12] Carroll C.: The Created Response Surface Technique for Optimizing Nonlinear-restrained Systems OR 9(2): 169-184, (1961)

[13] Censor Y., Zenios S.: The proximal minimization algorithm with d–functions," Journal of Optimization Theory and Applications, 73, 451–464, (1992)

[14] Chen C. and Mangasarian O. L.: Smoothing methods for convex inequalities and linear complementarity problems. Mathematical Programming 71, 51-69, (1995)

[15] Chen G. and Teboulle M.: Convergence analysis of a proximal–like minimization algorithm using Bregman Functions," SIAM J. Optimization, 3 (4), pp. 538–543. (1993)

[16] Courant R.: Variational methods for the solution of problems of equilibrium and vibrations. Bulletin of the American Mathematical Society, 49 1-23, (1943)

[17] Daube-Witherspoon M. and Muehllehner: An iterative space reconstruction algorithm suitable for volume ECT, IEEE Trans.Med Imaging M-5, 61-66, (1986)

[18] Dikin I.: Iterative Solutions of Linear and Quadratic Programming Problems," Soviet Mathematic Doklady, **8** 674–675, (1967)

[19] Eckstein J.: Nonlinear proximal point algorithms using Bregman functions with applications to convex programming," Mathematics of Operations Research, 18 (1), 202–226, (1993)

[20] Eggermont P.: Multiplicative iterative algorithm for convex programming. Linear Algebra and its Applications 130, 25-32, (1990)

[21] Ekeland I.: Legendre Duality in Nonconvex Optimization and Calculus of Variations.SIAM J. Control and Optimization. Vol. 16 No. 6 pp. 905-934 (1977)

[22] Fiacco A. Mc Cormick G.: Nonlinear programming, Sequential Unconstrained Minimization Techniques" SIAM, (1990)

[23] Frisch K.: The logarithmic Potential Method for convex programming. Memorandum of may 13 1955, University Institute of Economics, Oslo, (1955)

[24] Goldshtein E., Tretiakov N.: Modified Lagrangian functions. Moscow, (1989)

[25] Griva I., Polyak R.: Primal-dual nonlinear rescaling method with dynamic scaling parameter update. Math. Program. Ser. A 106, 237-259 (2006)

[26] Griva I., Polyak R.: Proximal Point Nonlinear Rescaling Method for Convex Optimization. Numerical Algebra, Control and Optimization, 1, N 3, 283-299, (2013)

[27] Guler O.: On the convergence of the proximal point algorithm for convex minimization. SIAM J. Control Optim., vol 29, pp 403–419, (1991)

[28] Hestenes M R.: Multipliers and gradient methods," JOTA, vol. 4, pp. 303–320, (1969)

[29] Hiriat-Urruty J. and Martinez-Legaz J.: New Formulas for the Legendre-Fenchel Transform. J. Math. Anal. Appl. 288 544-555 (2003)

[30] Ioffe A., Tichomirov V.: Duality of convex functions and extremum problems .Uspexi Mat.Nauk vol 23,n 6(144) ,51-116 (1968).

[31] Jensen D., Polyak R.: The convergence of a modify barrier method for convex programming IBM. Journal Res. Develop 38, 3 307-321, (1999)

[32] Kocvara M., Stingl M,.: PENNON. A code for convex nonlinear and semidefinite programming. Optimization methods and software, 18, 3, 317-333 (2003)

[33] Kocvara M., Stingl M.: Resent progress in the NLP-SDP code PENNON, Workshop "Optimization and Applications", Oberwalfach (2005)

[34] Kocvara M., Stingl M.: On the solution of large-scale SDP problems by the modified barrier method using iterative solver, Math. Program. series B, 109, 413-444 (2007)

[35] Martinet B.:Regularization d'inequations variationelles par approximations successive. Rev. Fr. Inf. Res. Ofer, V4, NR3 pp. 154-159, (1970)

[36] Martinet B.: Determination approachee d'un point fixe d'une application pseudo-contractante. C.R. Acad. Sci. Paris V274, N2 pp. 163-165, (1972)

[37] Melman A., Polyak R.: "The Newton modified barrier method for QP problems". Annals of OR ., vol. 62, pp. 465-519 (1996)

[38] Moreau J.: Proximite et dualite dans un espace Hilbertien. Bull. Soc. Math. France, V93 pp. 273-299, (1965)

[39] Matioli L., Gonzaga C.: A new family of penalties for Augmented Lagrangian methods," Numerical Linear Algebra with Applications, 15, 925–944, (2008)

[40] Motzkin T.: New Techniques for Linear Inequalities and Optimization. In project SCOOP, Symposium on Linear Inequalities and Programming, Planning Research Division, Director of Management Analysis Service, U.S. Air Force, Washington, D.C., no. 10, (1952)

[41] Nash S., Polyak R. and Sofer A.: "A numerical comparison of barrier and modified barrier method for large scale bound-constrained optimization". Large Scale Optimization, State of the Art. W. Hager, D. Hearn, P. Pardalos (Eds.). Kluwer Academic Publishers, pp. 319-338 (1994)

[42] Nesterov Yu., Nemirovsky A.: Interior Point Polynomial Algorithms in Convex Programming. SIAM, Philadelphia (1994)

[43] Nesterov Yu.:Introductory Lectures on Convex Optimization. Kluwer Academic Publishers, Norwell, MA (2004)

[44] Polyak B.: Iterative Methods Using Lagrange Multipliers for Solving Extremal Problems with Constraints of the Equation type. Comput. Math and Math Phys. 10, N 5, (1970)

[45] Polyak B.: Introduction to Optimization," Optimization Software, New York, NY, (1987)

[46] Polyak R.: Modified Barrier Functions (thory and methods). Math Programming 54, 177-222 (1992)

[47] Polyak R.: Log-Sigmoid Multipliers Method in Constrained Optimization. Annals of Operations Research 101, 427-460, (2001)

[48] Polyak R.: Nonlinear rescaling vs. smoothing technique in convex optimization, Math. Program. 92, 197–235, (2002)

[49] Polyak R.: Lagrangian Transformation and interior ellipsoid methods in Convex Optimization. JOTA V163, 3, (2015)

[50] Polyak R, Teboulle M.: Nonlinear rescaling and proximal–like methods in convex optimization. Mathematical Programming, 76, 265–284, (1997)

[51] Polyak B., Tret'yakov N.: The Method of Penalty Estimates for Conditional Extremum Problems. Comput. Math. and Math. Phys. 13, N 1, 42-58, (1973)

[52] Powell M. J. D.: A method for nonlinear constraints in minimization problems," in Fletcher (Ed.), Optimization, London Academic Press, pp. 283–298, (1969)

[53] Powell M.: Some convergence properties of the Modified Log Barrier Methods for Linear Programming. SIAM Journal on Optimization, vol 50 no. 4, 695-739, (1995)

[54] Ray A., Majumder S.: Derivation of some new distributions in statistical mechanics using maximum entropy approach. Yugoslav Journal of OR (24), NI, 145-155 (2014)

[55] Reich S. and Sabach S.: Two strong convergence theorems for a proximal method in reflexive Banach spaces. Numer. Funct. Anal. Optim. 31, 22-44, (2010)

[56] Rockafellar R. T.: A dual approach to solving nonlinear programming problems by unconstrained minimization," Math. Programming, 5, pp. 354–373, (1973)

[57] Rockafellar R. T.: Augmented Lagrangians and applications of the proximal points algorithms in convex programming," Math. Oper. Res., 1, 97–116, (1976)

[58] Rockafellar R. T.: Monotone operators and the proximal point algorithm," SIAM J. Control Optim, 14, pp. 877–898, (1976)

[59] Teboulle M.: Entropic proximal mappings with application to nonlinear programming. Mathematics of Operations Research 17, 670-690, (1992)

[60] Tikhonov A.N.: Solution of incorrectly formulated problems and the regularization method. Translated in Soviet Mathematics 4: 1035-1038, (1963)

[61] Tseng P., Bertsecas D.: On the convergence of the exponential multipliers method for convex programming. Math. Program, 60, 1-19 (1993)

[62] Vardi Y., Shepp L. and Kaufman L.: A statixtical model for position emission tomography. J. Amer. Statist. Assoc. 80, 8-38 (1985)

Department of Mathematics, The Technion - Israel Institute of Technology, 32000 Haifa, Israel

*E-mail address*: `rpolyak@techunix.technion.ac.il and rpolyak@gmu.edu`